\documentclass[11pt]{amsart}
\usepackage{latexsym}
\usepackage{amsmath}
\usepackage{amsfonts}
\usepackage{amssymb}
\usepackage{mathtools}
\usepackage[all]{xy}
\usepackage{pgf,tikz}
\usepackage{caption}
\usetikzlibrary{arrows}
\usetikzlibrary{decorations.pathmorphing}
\theoremstyle{plain}
\newcommand{\midarrow}{\tikz \draw[-triangle 90] (0,0) -- +(.1,0);}
\newcommand{\figref}[1]{\figurename~\ref{#1}}

\newcommand{\cleqn}{\setcounter{equation}{0}}
\newcommand{\clth}{\setcounter{theorem}{0}}
\newcommand {\sectionnew}[1]{\section{#1}\cleqn\clth}

\newtheorem{theorem}{Theorem}[section]
\newtheorem{lemma}[theorem]{Lemma}
\newtheorem{definition-theorem}[theorem]{Definition-Theorem}
\newtheorem{proposition}[theorem]{Proposition}
\newtheorem{corollary}[theorem]{Corollary}
\newtheorem{definition}[theorem]{Definition}
\newtheorem{example}[theorem]{Example}
\newtheorem{remark}[theorem]{Remark}
\newtheorem{conjecture}[theorem]{Conjecture}

\newcommand \bth[1] { \begin{theorem}\label{t#1} }
\newcommand \ble[1] { \begin{lemma}\label{l#1} }

\newcommand \bpr[1] { \begin{proposition}\label{p#1} }
\newcommand \bco[1] { \begin{corollary}\label{c#1} }
\newcommand \bde[1] { \begin{definition}\label{d#1}\rm }
\newcommand \bex[1] { \begin{example}\label{e#1}\rm }
\newcommand \bre[1] { \begin{remark}\label{r#1}\rm }
\newcommand \bcj[1] { \begin{conjecture}\label{j#1}\rm }

\renewcommand {\eth} { \end{theorem} }
\newcommand {\ele} { \end{lemma} }

\newcommand {\epr} { \end{proposition} }
\newcommand {\eco} { \end{corollary} }
\newcommand {\ede} { \end{definition} }
\newcommand {\eex} { \end{example} }
\newcommand {\ere} { \end{remark} }
\newcommand {\ecj} { \end{conjecture} }

\newcommand {\enota} { \end{notation} }
\newcommand \thref[1]{Theorem \ref{t#1}}
\newcommand \leref[1]{Lemma \ref{l#1}}
\newcommand \prref[1]{Proposition \ref{p#1}}
\newcommand \coref[1]{Corollary \ref{c#1}}

\newcommand \deref[1]{Definition \ref{d#1}}

\newcommand \reref[1]{Remark \ref{r#1}}
\newcommand \lb[1]{\label{#1}}

\def \Zset {{\mathbb Z}}

\def \Qset {{\mathbb Q}}

\def \AA  {{\mathcal{A}}}           

\def \CC {{\mathcal{C}}}

\def \QQ {{\mathcal{Q}}}
\def \FF {{\mathcal{F}}}

\def \UU {{\mathcal{U}}}

\def \TT {{\mathcal{T}}}

\def \De {\Delta}   
\def \de {\delta}
\def \al {\alpha}
\def \be {\beta}

\def \ga {\gamma}
\def \de {\delta}

\def \sig {\sigma}

\def \sig{\sigma}

\def \mt  {\mapsto}
\def \lra {\longrightarrow}

\def \hra {\hookrightarrow}

\def \ol {\overline}
\def \wt {\widetilde}

\DeclareMathOperator \Aut { {\mathrm{Aut}} }

\DeclareMathOperator \Homeo { {\mathrm{Homeo}}}
\DeclareMathOperator \MCG { {\mathcal{MCG}}}
\DeclareMathOperator \Mod {{\mathcal{M}}od}

\DeclareMathOperator \opp { {\mathrm{op}}}

\DeclareMathOperator \Iso { {\mathrm{Iso}} }

\begin{document}
\title[Recovering the topology of surfaces from  cluster algebras]
{Recovering the topology of surfaces \\ from  cluster algebras}
\author[Eric Bucher]{Eric Bucher}
\address{
Department of Mathematics \\
Louisiana State University \\
Baton Rouge, LA 70803
U.S.A.
}
\email{ebuche2@math.lsu.edu}
\author[Milen Yakimov]{Milen Yakimov}
\email{yakimov@math.lsu.edu}
\thanks{The research of E.B. has been supported by a GAANN fellowship and a VIGRE fellowship through the NSF grant DMS-0739382 
and that of M.Y. by the NSF grants DMS-1303038 and DMS-1601862.}
\begin{abstract} We present an effective method for recovering the topology of a bordered oriented surface with marked points from its cluster algebra. 
The information is extracted from the maximal triangulations of the surface, those that have exchange quivers with maximal number of arrows in the mutation 
class. The method gives new proofs of the automorphism and isomorphism problems for the surface cluster algebras as well as the 
uniqueness of the Fomin--Shapiro--Thurston block decompositions of the exchange quivers of the surface cluster algebras. 
The previous proofs of these results followed a different approach based on Gu's direct proof of the last result. 
The method also explains the exceptions to these results due to pathological problems with the maximal 
triangulations of several surfaces.
\end{abstract}
\definecolor{zzttqq}{rgb}{0.6,0.2,0}
\maketitle
\sectionnew{Introduction}
\lb{Intro}
\subsection{Cluster algebras} 
\label{1.1}
Cluster algebras form a large, axiomatically defined class of algebras, introduced by Fomin and Zelevinsky 
\cite{FZ1,FZ2}. They play a key role in many settings in Lie and representation theory, combinatorics, mathematical physics, algebraic and Poisson geometry,
see \cite{F,Le,L,M,W} for details. 

A cluster algebra $\AA$ (without coefficients) is defined by recursively mutating finite sets of algebraically independent elements $(x_1, \ldots, x_n)$ 
of $\AA$, called clusters. The mutations at each step are determined by a quiver $Q$ with $n$ vertices, called an exchange quiver, 
and more generally by an integer matrix. Cluster algebras with coefficients  involve additional mutation in semifields \cite{FZ1}.
The pairs $((x_1, \ldots, x_n), Q)$ are called seeds of $\AA$.

An important subclass of cluster algebras establishes a link with Teichm\"uller theory. To every oriented surface 
$S$ (with or without boundary) and a set of marked points $M$ on it, Gekhtman--Shapiro--Vainshtein \cite{GSV}, Fock--Goncharov \cite{FG},
and Fomin--Shapiro--Thurston \cite{FST,FT}, associated a cluster algebra $\AA(S,M)$. A detailed study of their 
cluster variables and seeds was carried out in \cite{FST,FT} in terms of tagged arcs and triangulations. 
We will refer to the pair $(S,M)$ as a {\em{bordered surface with marked points}}.

Felikson, Tumarkin and Shapiro \cite{FSTu} proved that the surface cluster algebras and several exceptional algebras exhaust 
all skewsymmetric cluster algebras of finite mutation type. 
Surface cluster algebras have been studied from many different perspectives and have been used in various applications.
The cluster expansions in them were described in \cite{MSW1} and bases of the cluster algebras were constructed in \cite{MSW2,CLR}.
Their additive categorifications and the representation type of the involved algebras were studied in \cite{La,GLS,L}.   
The local acyclicity of these cluster algebras was obtained in \cite{Mu} and the relation between cluster and upper cluster 
algebras was described in \cite{Mu,CLR}.
The surface cluster algebras were applied to the study of BPS states and spaces of stability conditions in \cite{ACC,BS}, 
and to the problem of constructing presentations of 
Coxeter and braid groups \cite{BM,FeT,GM}. The problem of existence of maximal green sequences was treated in \cite{ACC,BMi,L2}.
Also, the quantum counterparts were related to skein algebras in \cite{Mu2}. 
\subsection{The automorphism and isomorphism problems, uniqueness of the Fomin--Shapiro-Thurston block decompositions} 
\label{1.2}
Given a class of algebras $\CC$, two main problems are: {\em{the Isomorphism Problem}} (for $A, B \in \CC$ determine when $A \cong B$) and
{\em{the Automorphism Problem}} (for $A \in \CC$, describe $\Aut A$). 

Fomin and Zelevinsky \cite{FZ2} defined a {\em{strong cluster isomorphism}} between the cluster algebras $\AA_1$ and $\AA_2$ to be any 
algebra isomorphism $\phi \colon \AA_1 \stackrel{\cong}\lra \AA_2$ that has the property that there is a seed $((x_1, \ldots, x_n), Q)$ of 
$\AA_1$ such that $((\phi(x_1), \ldots, \phi(x_n)), Q)$ is a seed of $\AA_2$. This immediately implies that any seed of $\AA_1$ will 
have this property. Denote the set of such isomorphisms by $\Iso^+(\AA_1, \AA_2)$.

Assem, Schiffler and Shramchenko \cite{ASS} defined a {\em{cluster isomorphism}} between the cluster algebras $\AA_1$ and $\AA_2$ to be an 
algebra isomorphism $\phi \colon \AA_1 \stackrel{\cong}\lra \AA_2$ that has the property that there is a cluster $(x_1, \ldots, x_n)$ of 
$\AA_1$ such that $(\phi(x_1), \ldots, \phi(x_n))$ is a cluster of $\AA_2$ and $\phi$ commutes with mutation in the sense that $\phi(\mu_k (x_k)) = \mu_k (\phi(x_k))$ for all $k$, where $\mu_k$ denotes the mutation in the direction $k$.  
When $Q$ is connected, this is equivalent to saying that 
$\phi$ is an algebra isomorphism with the property that there is a seed $((x_1, \ldots, x_n), Q)$ of 
$\AA_1$ such that either $((\phi(x_1), \ldots, \phi(x_n)), Q)$ or $((\phi(x_1), \ldots, \phi(x_n)), Q^{\opp})$ is a seed of $\AA_2$, where 
$Q^{\opp}$ denotes the opposite quiver. Once again, it is an easy corollary of the definition of mutation that all clusters and seeds of $\AA_1$ 
will posses the stated properties with respect to the map $\phi$. Denote the set of these isomorphisms by $\Iso(\AA_1, \AA_2)$. Let
\[
\Aut^+ (\AA) := \Iso^+ (\AA, \AA) \quad \mbox{and} \quad \Aut(\AA) := \Iso (\AA, \AA). 
\]

Fomin, Shapiro and Thurston proved that, the exchange quiver $Q_{\TT}$ of every triangulation of a surface 
$(S,M)$ can be obtained by gluing blocks and then canceling pairs of opposite arrows. Each such presentation of $Q_{\TT}$ is 
called a {\em{block decomposition}}. With two exceptions, the block decomposition of $Q_{\TT}$ comes from a decomposition of 
the triangulation $\TT$ into {\em{puzzle pieces}}, see for details \S \ref{4.2}. 

The problem of {\em{recovering the topology of a surface from its cluster algebra structure}} was raised by 
Fomin, Shapiro and Thurston \cite[Sect. 14]{FST}. They proved \cite[Prop. 14.1]{FST} that, if
$\TT_1$ and $\TT_2$ are triangulations of bordered oriented surfaces with marked points $(S_1,M_1)$ 
and $(S_2, M_2)$, such that $Q_{\TT_1} \cong Q_{\TT_2}$ and $Q_{\TT_1}$ has a unique block decomposition, then 
$(S_1, M_1)$ and $(S_2, M_2)$ are homeomorphic. This relates the isomorphism and automorphism problems for the class of 
surface cluster algebras to {\em{the problem for uniqueness of block decompositions of exchange quivers}}. 

Weiwen Gu proved \cite{G1,G1b,G2} that, with six exceptions, the exchange quivers of surface cluster algebras have unique
block decompositions. This proof is based on a direct and long argument for the exchange quivers of all possible triangulations 
of surfaces, not special triangulations. The isomorphism problem for surface cluster algebras is solved from here by applying 
\cite[Prop. 14.1]{FST}. Bridgeland and Smith \cite{BS} solved the automorphism 
problem, based on Gu's result. 

In this paper, we present an effective method for recovering the topology of a surface from its cluster algebra. 
The main idea is to consider only ``maximal triangulations'' -- those whose exchange quivers have maximal number of arrows, 
and then to derive results for all triangulations from the maximal ones.  
This gives conceptual new proofs of the uniqueness of block decompositions, the isomorphism and automorphism 
problems for surface cluster algebras. 
\subsection{Maximal triangulations}
\label{1.3}
We call a triangle in a triangulation {\em{a face}}, {\em{a wedge}}, and {\em{a cap}} if it has 0, 1, and 2 of its sides on the boundary of the surface,
respectively. 

The number of edges $e(Q_\TT)$ of a triangulation of a bordered oriented surface with marked points $(S,M)$ has an easy upper bound in 
terms of the combinatorial data (see \coref{max-triang} for details). We say that a triangulation $\TT$ is {\em{maximal}} if this upper bound is attained.
The exchange quiver $Q_\TT$ of such a triangulation has maximal number of edges in its mutation class. 

It is straightforward to show that maximal triangulations are precisely the triangulations $\TT$ that satisfy the following three conditions:
\begin{enumerate}
\item[(a)] $\TT$ does not contain pairs of triangles glued along 2 edges in a way that 
a pair of arrows cancels out in the computation of $Q_\TT$  (we call those ``negatively double-glued pairs of triangles''), 
\item[(b)] $\TT$ does not have self-folded triangles, and 
\item[(c)] $\TT$ has maximal possible number of caps. 
\end{enumerate}
Any surface that is different from the once-punctured digon and the twice-punctured monogon has a maximal triangulation

A triangulation $\TT$ of $(S,M)$ will be called {\em{connected}}, if one can get from any face of $\TT$ to any other face by a path 
crossing only faces and no marked points. 

Our key result is that one can reconstruct $(S,M)$ from the exchange quiver $Q_\TT$ of any 
{\em{maximal connected triangulation $\TT$ with at least two faces}}. Only very few surfaces do 
not have such triangulations (see \thref{exist}) -- this is one of the features that make our method effective. 
\medskip
\\
\noindent
{\bf{Theorem A.}} {\em{Let $(S_1, M_1)$ and $(S_2, M_2)$ be two bordered surfaces with marked points which are different 
from the 4-punctured sphere, and the twice-punctured monogon and digon.

Assume that $\TT_1$ and $\TT_2$ are triangulations of $(S_1, M_1)$ and $(S_2, M_2)$ such that $\TT_1$ is  
a connected maximal triangulation with at least 2 faces, and that
$\phi \colon Q_{\TT_1} \stackrel{\cong}\lra Q_{\TT_2}$ is an isomorphism of quivers. Then there exist an orientation-preserving 
homeomorphism $g \colon (S_1, M_1) \stackrel{\cong}\lra (S_2, M_2)$ which induces $\phi$}}.
\medskip
\\
\noindent 
We prove the theorem in the following steps:

(I) Since $Q_{\TT_2}$ has maximal number of arrows in its mutation class, $\TT_2$ would have to be maximal too.

(II) We show that, if two faces $\De'$ and $\De''$ of $\TT_1$ share at least one 
common side, then the isomorphism $\phi \colon Q_{\TT_1} \stackrel{\cong}\lra Q_{\TT_2}$ takes the vertices of $Q_{\TT_1}$, corresponding to 
the arcs of $\De'$ and $\De''$, to vertices in $Q_{\TT_2}$ that come from the arcs of a pair of faces in $\TT_2$ glued in the same way as $\De'$ and $\De''$. 
We recursively apply this to construct a homeomorphism $g$ from the union of all faces of $\TT_1$ (without the marked points $M_1$) 
to the union of faces of $\TT_2$ (without the points $M_2$).

(III) We extend $g$ to $S_1$ by recursively adding the wedges and caps of $\TT_1$.

A more general result for tagged traingulations $\TT_i$ is proved in \thref{homeo}.
The notion of a connected maximal triangulation with at least 2 faces is related to the notion 
of a gentle triangulation of Geiss, Labradini-Fragoso, and Schr\"oer \cite{GLS}, but the two notions are different, see \S \ref{4.3} for details. 
The latter notion is more restrictive and leads to a larger number of exceptions, including three infinite families,
\cite[Theorem 7.8]{GLS}.
\subsection{Consequences} 
\label{1.4}
Define the {\em{mapping class group}} of a bordered surface with marked points $(S,M)$ by 
\[
\MCG(S,M):= \Homeo^+(S,M) / \Homeo^+_0(S,M)
\]
where $\Homeo^+(S,M)$ is the group of orientation-preserving homeomorphisms $g \colon S \stackrel{\cong}\lra S$ 
such that $g(M)= M$, and $\Homeo^+_0(S,M)$ is its subgroup of homeomorphisms homotopic to the identity. 
The group $\MCG(S,M)$ acts on the set of punctures $P$ of $(S,M)$. The 
{\em{tagged mapping class group}} of $(S,M)$ is defined \cite{ASS,BS} to be the 
semidirect product 
\[
\MCG_{\bowtie} (S,M) := \MCG(S,M) \ltimes \Zset_2^P.
\]
Theorem A (or more precisely, its tagged generalization in \thref{homeo}) implies that, with few exceptions, the elements of $\Iso^+(\AA(S_1, M_1), \AA(S_2, M_2))$ 
are induced by orientation-preserving homeomorphisms $g \colon (S_1, M_1) \stackrel{\cong}\to (S_2, M_2)$. 
This immediately gives conceptual new proofs of the following results:
\medskip
\\
\noindent
{\bf{Theorem B.}} (Gu, Fomin--Shapiro--Thurston) \cite{G1, FST} {\em{Let $(S_1,M_1)$ and $(S_2, M_2)$ be two bordered surfaces with marked points. Then, 
$\Iso^+(\AA(S_1, M_1), \AA(S_2, M_2)) \neq \varnothing$ if and only if $\Iso(\AA(S_1, M_1), \AA(S_2, M_2))  \neq \varnothing$.

Furthermore, $\Iso^+(\AA(S_1, M_1), \AA(S_2, M_2)) \neq \varnothing$ if an only if 
there exists a homeomorphism $g \colon (S_1, M_1) \stackrel{\cong}\lra (S_2, M_2)$
or the pair of surfaces is one of the following:

(a) the unpunctured hexagon and once-punctured triangle;

(b) the the twice-punctured monogon and the annulus with $(2,2)$ marked points on the boundary.}}
\medskip
\\
\noindent
{\bf{Theorem C.}} (Gu, Bridgeland--Smith) \cite{G1,BS} {\em{Let $(S,M)$ be a bordered surface with marked points which is different from 
the 4-punctured sphere, the once-punctured 4-gon and the twice-punctured digon. 

(a) If $(S,M)$ is not a once-punctured closed surface, then
\[
\Aut^+ \AA(S,M) \cong \MCG_{\bowtie}(S,M). 
\]

(b) If $(S,M)$ is a once-punctured closed surface, then
\[
\Aut^+ \AA(S,M) \cong \MCG(S,M). 
\]
}}
\medskip

If an exchange quiver of a surface cluster algebra has two block decompositions, then there exist two surfaces $(S_1,M_1)$,
$(S_2,M_2)$ and a strong cluster isomorphism $\phi \in \Iso^+(\AA(S_1, M_1), \AA(S_2, M_2))$. 
With a few exceptions we can choose a connected maximal triangulation $\TT_1$ of $(S_1,M_1)$ 
with at least 2 faces and apply Theorem A to it and $\phi$ (or more precisely, the tagged analog in 
\thref{homeo}). This implies that $\phi$ is induced by a homeomorphism $g \colon (S_1, M_1) \to (S_2,M_2)$ 
and thus the two block decompositions that we started with are identical. This gives a new proof of the following theorem:
\medskip
\\
\noindent
{\bf{Theorem D.}} (Gu) \cite{G1} {\em{Let $(S,M)$ be a bordered surface with marked points which is different from one of the following: 

The 4-punctured sphere, the twice-punctured digon, the once-punctured 2-, 3- and 4-gons,
the annulus with $(2,2)$ marked points on the boundary and the unpunctured hexagon.

For every triangulation $\TT$ of $(S,M)$, its exchange quiver $Q_\TT$ has a unique block decomposition.
}}
\medskip

Our approach to Theorems B--D via maximal triangulations sheds light on the intrinsic reason for the 
uniqueness of block decompositions and explains the exceptions in those theorems. They are due to two
facts: the lack of maximal triangulations for a few surfaces and a pathological problem coming from the impossibility to match 
pairs of glued faces in the three cases described in \prref{match2}.

In \cite{Fr} Fraser defined and studied a notion of quasi-homomorphisms between cluster algebras with arbitrary coefficients 
from semifields. He proved a theorem classifying the quasi-automorphism groups of surface cluster algebras with arbitrary coefficients 
on the basis of the coefficient free result from Theorem C. 
For cluster algebras with coefficients of geometric type, cluster morphisms were also defined and studied in \cite{ADS,LS}.

This paper is organized as follows: Sections \ref{backgr-scl} and \ref{backgr-auto} contain background material on surface cluster algebras and cluster automorphism
groups. Section \ref{match-fa} defines and describes the main properties of maximal triangulations, including
Step I and the first part of Step II of the proof of Theorem A; the recursive Steps II-III of the proof are  carried out in Section \ref{homeo}.
Theorems B/D and C are proved in Sections \ref{cl-iso} and \ref{cl-auto}, respectively. 
\medskip
\\
\noindent
{\bf Acknowledgements.} We are grateful to Sergey Fomin and Misha Shapiro for their very helpful suggestions and comments 
on the first version of the paper. We are indebted to the referee for pointing out inaccuracies and making many suggestions 
which improved the paper.
M. Y. would like to thank Newcastle University and the Max Planck Institute for Mathematics in Bonn
for the warm hospitality during visits in the Fall of 2015.
\sectionnew{Background on surface cluster algebras}
\lb{backgr-scl}
This section contains background material on cluster algebras and the construction of cluster algebras from bordered 
surfaces with marked points. 
\lb{CA-back}
\subsection{Cluster Algebras}
\label{2.1} Cluster algebras of skewsymmetric (geometric type without frozen variables) 
are defined starting from a quiver (a directed graph) without loops and 2-cycles. A quiver 
is represented as a quadruple $Q= (Q_0, Q_1, s, t)$ where $Q_0$ is the set of vertices, $Q_1$ is the set of arrows, 
and $s, t \colon Q_1 \to Q_0$ are the source and target maps for the arrows. Let 
\[
n:= |Q_0| \quad \mbox{and} \quad e(Q) := |Q_1|.
\]
Denote the 
set of vertices $Q_0 = \{1, \ldots, n\}$. Consider the ambient field $\FF:= \Qset(x_1, \ldots, x_n)$
and set ${\bf{x}}:= (x_1, \ldots, x_n)$. The cluster algebra $\AA({\bf{x}}, Q)$ without frozen variables is a $\Zset$-subalgebra 
of $\FF$ generated by a (generally infinite) set of cluster variables obtained by an iterative procedure of mutation. 

For $1 \leq k \leq n$, the {\em{mutated quiver}} $\mu_k(Q)$ is obtained from $Q$ by performing the following operations:
\begin{enumerate}
\item for each subquiver $i \to k \to j$, add a new arrow $i \to j$;
\item reverse all allows with source or target $k$;
\item remove the arrows in a maximal set of pairwise disjoint 2-cycles.
\end{enumerate}
The pair $({\bf{x}}, Q)$ is called an {\em{initial seed}}. Its {\em{mutations}} are defined by 
\[
\mu_k({\bf{x}}, Q):= ({\bf{x}}', \mu_k(Q)), \quad {\bf{x}}' := ({\bf{x}} \backslash \{x_k\}) \cup \{x'_k\}
\]
where
\[
x'_k := \frac{1}{x_k} \left( \prod_{\al \in Q_1, s(\al)=k} x_{t(\al)} + \prod_{\al \in Q_1, t(\al)=k} x_{s(\al)} \right).
\]
All pairs obtained by successive mutations from $(\bf{x}, Q)$ are called {\em{labeled seeds}}. They have the property that each 
cluster $\wt{\bf{x}}$ is a trancendence
basis of $\FF$ over $\Qset$. The elements of $\wt{\bf{x}}$ are called {\em{cluster variables}}. The cluster algebra $\AA({\bf{x}}, Q)$
is the $\Zset$-subalgebra of $\FF$ generated by all cluster variables. The number $n = |\wt{\bf{x}}|$ is called rank of the cluster 
algebra. 

Every seed $(\wt{\bf{x}}, \wt{Q})$ of the cluster algebra $\AA({\bf{x}}, Q)$ is uniquely determined by its cluster $\wt{\bf{x}}$, \cite{GSV2}. For simplicity 
of the notation, the vertices of the quiver $\wt{Q}$ will be indexed by the cluster variables in $\wt{\bf{x}}$.
\subsection{Surface cluster algebras}
\label{2.2} Let $S$ be a {\em{connected oriented Riemann surface}} with or without boundary, and $M \subset S$ be a finite set of {\em{marked points}} such that 
each connected component of the boundary $\partial S$ contains a marked point. The marked points in the interior of $S$ will be 
called {\em{punctures}}. 

\bde{bor} A {\em{bordered surface with marked points}} is a pair $(S,M)$ as above that is not one of the following: 

{\em{A sphere with 1, 2 or 3 marked points, 
an unpunctured disk with 1,2 or 3 marked points on the boundary, or a once-punctured disk with 1 marked 
point on the boundary.}}
\ede

\bde{arc} (a) An {\em{arc}} in $(S,M)$ is the isotopy class of a curve in $S$ connecting 2 marked points such that
\begin{itemize} 
\item the curve does not have self-intersections, except possibly coinciding end points;

\item the interior of the curve is disjoint from $M$ and $\partial S$;

\item the curve does not cut out an unpuctured monogon or an unpunctured digon.
\end{itemize}
The set of arcs of $(S,M)$ will be denoted by $A(S,M)$.

(b) Two arcs are called {\em{compatible}} if the their classes contain curves which do not intersect, except possibly at the end points.

(c) A maximal collection of distinct compatible arcs is called an {\em{ideal triangulation}} (or simply a {\em{triangulation}}).
\ede
{\em{We will use the following:}}
\medskip
\\
\noindent
{\bf{Convention for triangulations.}} For each triangulation $\TT$ of a bordered surface with marked points $(S,M)$, 
one can choose representatives of the arcs of $\TT$ that do not intersect each other, 
except possibly at the end points of the arcs.
{\em{We fix once and for all such a presentation of each triangulation of $(S,M)$}}. It will be used for inductive 
constructions of homeomorphisms. 

\bde{triangles} 
(a) By an {\em{arc}} of the triangulation $\TT$ we will mean the particular {\em{closed curve}} from the presentation of $\TT$. Denote 
the set of arcs of $\TT$ by $A(\TT)$.

(b) Denote by $BA(S,M)$ the closures of the connected components of $\partial S \backslash M$. The elements of $BA(S,M)$ will be called 
{\em{boundary arcs}} of $(S,M)$.

(c) By a triangle of $\TT$ we will mean the closure of a connected component of the complement to the set of arcs of $\TT$ in $S$. 
Denote the set of triangles of a triangulation $\TT$ by $T(\TT)$. 

(d) A triangle will be called a {\em{face}}, a {\em{wedge}} or a {\em{cap}} if its boundary contains exactly 0, 1, or 2 boundary arcs, respectively.

(e) A {\em{self-folded triangle}} is a face with 2 sides (arcs) that coincide. The repeated arc is called a {\em{radius}} and the other one is called a {\em{loop}}.  
\ede

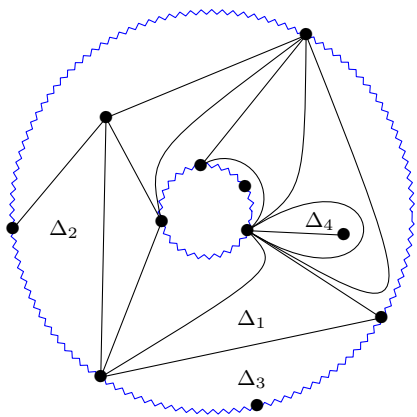
\begin{figure}
\begin{center}
\begin{tikzpicture}[line cap=round,line join=round,>=triangle 45,x=1.0cm,y=1.0cm]
\clip(-4.3,-0.96) rectangle (5.66,6.3);
\draw(0.47,2.73)[color=blue, line join=round, decorate, decoration={zigzag,
    segment length=4,
    amplitude=.9,post=lineto,
    post length=2pt
}] circle (2.66cm);
\draw(0.39,2.74)[color=blue, line join=round, decorate, decoration={zigzag,
    segment length=4,
    amplitude=.9,post=lineto,
    post length=2pt
}] circle (0.61cm);
\draw [shift={(0.63,2.91)}] plot[domain=-0.95:2.2,variable=\t]({1*0.53*cos(\t r)+0*0.53*sin(\t r)},{0*0.53*cos(\t r)+1*0.53*sin(\t r)});
\draw (-1.01,0.53)-- (-0.94,3.98);
\draw (-0.94,3.98)-- (-2.18,2.5);
\draw (-0.94,3.98)-- (-0.2,2.6);
\draw (0.32,3.34)-- (1.72,5.08);
\draw (0.94,2.48)-- (2.72,1.32);
\draw (-0.94,3.98)-- (1.72,5.08);
\draw (-1.01,0.53)-- (-0.2,2.6);
\draw (-1.01,0.53)-- (2.72,1.32);
\draw (0.94,2.48)-- (2.22,2.42);
\draw (0.94,2.48).. controls (1.7,3) ..(1.72,5.08);
\draw (-0.2,2.6).. controls (-.4,3.5) ..(1.72,5.08);
\draw (-1.01,0.53).. controls (1.4,2) ..(0.94,2.48);
\draw (0.94,2.48).. controls (3,3.8) and (3,1.2)..(0.94,2.48);
\draw (0.94,2.48).. controls (3.3,1.1) ..(1.72,5.08);
\begin{scriptsize}
\fill  (1.72,5.08) node {\large{$\bullet$}}; 
\fill (-2.18,2.5)node {\large{$\bullet$}}; 
\fill  (2.72,1.32) node {\large{$\bullet$}}; 
\fill (0.32,3.34)node {\large{$\bullet$}}; 
\fill (-0.2,2.6)node {\large{$\bullet$}}; 
\fill (0.94,2.48) node {\large{$\bullet$}}; 
\fill (0.91,3.06) node {\large{$\bullet$}}; 
\fill (-0.94,3.98) node {\large{$\bullet$}}; 
\fill (-1.01,0.53) node {\large{$\bullet$}}; 
\fill  (1.07,0.15) node {\large{$\bullet$}}; 
\fill (2.22,2.42) node {\large{$\bullet$}}; 
\draw (1,0.5) node{$\De_3$};
\draw (1,1.3) node{$\De_1$};
\draw (-1.5,2.5) node{$\De_2$};
\draw (1.9,2.6) node{$\De_4$};
\end{scriptsize}
\end{tikzpicture}
\end{center}
\caption{A triangulation
\label{extriang}}
\end{figure}
\figref{extriang} illustrates the different notions. Here and below solid black lines indicate arcs and wavy blue lines indicate boundary arcs. 
The triangles $\De_1$, $\De_2$, $\De_3$ and $\De_4$ are a face, wedge, cap and self-folded triangle, respectively.

Denote by $g$ the genus of the original Riemann surface, by $p$ the number of punctures, by $b$ the number of boundary components, 
and by $h_1, \ldots, h_b$ the number of marked points on the different boundary components. Set 
\[
h := h_1 + \cdots + h_b = |BA(S,M)|.
\]
The number of arcs of every triangulation $\TT$ of $(S,M)$ does not depend on the choice of triangulation, \cite{FST}:
\begin{equation}
\label{rank}
n = 6(g - 1) + 3 b + 3 p + h. 
\end{equation}
It is called the {\em{rank}} of $(S,M)$. When we need to specify the underlying $(S,M)$ in the above notations, we will 
write $n(S,M)$, etc. Denote 
\begin{equation}
\label{cSM}
c(S,M) := \sum_{i=1}^b \left\lfloor {\frac{h_i}{2}} \right\rfloor.
\end{equation}

The number of faces, wedges, and caps of a triangulation $\TT$ will be denoted by $f(\TT)$, $w(\TT)$, and $c(\TT)$, respectively.
We have, 
\begin{align}
3f(\TT) + 2 w(\TT) + c(\TT) &= 2 n.
\label{fwc1}
\\
w(\TT) + 2 c (\TT) &= h.
\label{fwc2}
\end{align}
The first equality follows from the fact that faces, wedges, and caps have 3, 2, and 1 arcs, respectively, 
and the fact that each arc is on the boundary of 2 triangles (both times counted with multiplicities). 
The second equality follows by counting the boundary arcs in wedges and caps, which equal to 1 and 2, 
respectively, and using that $h = |BA(S,M)|$. Finally, we also have 
an inequality on the number of caps of each triangulation 
\begin{equation}
\label{caps}
c(\TT) \leq c(S,M)
\end{equation}
since caps do not overlap.

To each triangulation $\TT$ of $(S,M)$, one associates an {\em{exchange quiver}} $Q_\TT$ whose vertices are indexed by the arcs $A(\TT)$.
{\em{To distinguish the arcs of $\TT$ from the vertices of $Q_\TT$, we will denote by $[\al]$ the vertex of $Q_\TT$ corresponding to the arc
$\al$.}} The edges of $Q_\TT$ are determined as follows:
\begin{enumerate}
\item For every triangle $\De \in T(\TT)$ that is not self-folded 
add an arrow $[\al] \to [\be]$ in each of the following cases:
\begin{itemize}
\item $\al$ and $\be$ are sides of $\De$, and $\be$ follows $\al$ in the clockwise order;
\item $\be$ is a radius of a self-folded triangle with a loop $\ga$, and $\al$ and $\ga$ are 
sides of $\De$ such that $\ga$ follows $\al$ in the clockwise order;
\item $\al$ is a radius of a self-folded triangle with a loop $\ga$, and $\be$ and $\ga$ are 
sides of $\De$ such that $\be$ follows $\ga$ in the clockwise order.
\end{itemize}
\item Remove the arrows in a maximal set of pairwise disjoint 2-cycles.
\end{enumerate}
The {\em{cluster algebra}}, $\AA(S,M)$ associated to the bordered surface with marked points $(S,M)$, 
is defined \cite{FG,FST,GSV} to be $\AA({\bf{x}},Q_\TT)$ for every triangulation $\TT$ of $(S,M)$. 
Each arc $\al \in \TT$ gives rise to a cluster variable $x_\al \in \AA(S,M)$ and each triangulation 
$\TT$ of $(S,M)$ gives rise to the seed $( \{x_\al \mid \al \in A(\TT) \}, Q_\TT)$ of $\AA(S,M)$. 
The process of mutation is represented by flips of arcs.
\subsection{Tagged triangulations}
\label{2.3} The seeds of the cluster algebras $\AA(S,M)$ are classified \cite{FST,FT} in terms of tagged triangulations of $(S,M)$
which generalize ordinary triangulations.  

A {\em{tagged arc}} is an arc on $(S,M)$ whose ends are marked (tagged) in 2 possible ways, {\em{plain}} or {\em{notched}}, so that the 
following conditions are satisfied:
\begin{itemize}
\item the arc does not cut out a once-punctured monogon; 
\item an endpoint lying on the boundary is tagged plain; 
\item if the arc is a loop, its endpoints are tagged in the same way.
\end{itemize}
The set of tagged arcs of $(S,M)$ will be denoted by $A_{\bowtie}(S,M)$.

Two tagged arcs $\al$ and $\be$ are called {\em{compatible}} if the plain arcs $\ol{\al}$ and $\ol{\be}$, obtained from $\al$ and $\be$ 
by forgetting the taggings, are compatible and satisfy the following:
\begin{itemize}
\item if $\ol{\al} = \ol{\be}$, then at least one end of $\al$ and $\be$ is tagged in the same way;
\item if $\ol{\al} \neq \ol{\be}$, but $\al$ and $\be$ have a common end point, then their taggings at this point
are the same.
\end{itemize}
A {\em{tagged triangulation}} of $(S,M)$ is a maximal collection of distinct pairwise compatible tagged arcs. Each 
tagged triangulation $\TT$ gives rise to an ordinary triangulation $\TT^\circ$ in the following way. The {\em{signature}} of a 
puncture $y$, with respect to a tagged triangulation $\TT$, is defined by
\[
\de_\TT(y) 
=
\begin{cases}
1, & \mbox{if all tagged arcs of $\TT$, containing $y$, are tagged plain at $y$}
\\
-1, & \mbox{if all tagged arcs of $\TT$, containing $y$, are tagged notched at $y$}
\\
0, & \mbox{otherwise}.
\end{cases}
\] 
The definition of tagged triangulation easily implies that in the third case there are precisely 2 arcs of $\TT$ containing $y$, $\al$ and $\be$, 
such that $\ol{\al}= \ol{\be}$ and the taggings of $\al$ and $\be$ at $y$ are different, while at the other end are the same. 
To each tagged triangulation $\TT$, one associates and ordinary triangulation $\TT^\circ$ by performing the two operations:
\begin{itemize}
\item replace all notched ends of arcs at the punctures with $\de_\TT(y)=-1$ by plain ones;
\item for each puncture $y$ with $\de_\TT(y) = 0$, we will have two arcs $\al$ and $\be$ containing $y$ which will satisfy $\ol{\al}= \ol{\be}$ 
and have different taggings at $y$ and the same taggings at the other endpoint; replace the arc $\be$ notched at $y$ with a loop based at the 
other end point of $\be$ and closely wrapping around $\be$.
\end{itemize}
The set of tagged arcs of $\TT$ will be denoted by $A_{\bowtie} (\TT)$. There is an obvious bijection
$A_{\bowtie}(\TT) \stackrel{\cong}\lra A(\TT^\circ)$. The vertices of the exchange quiver $\QQ_\TT$ will be indexed by $A_{\bowtie} (\TT)$; 
the vertex corresponding to $\al \in  A_{\bowtie} (\TT)$ will be denoted by $[\al]$. The edge set of $\QQ_\TT$ is defined by 
\[
Q_\TT:= Q_{\TT^\circ}
\]
in the above bijection.

We have an embedding $\tau \colon A(S,M) \hra A_{\bowtie}(S,M)$. The map sends every arc that is cutting a once-punctured monogon
to the radius of the corresponding self-folded triangle notched at the puncture of the monogon, and is the identity otherwise. This way, each ordinary triangulation 
$\TT$ gives rise to a tagged one $\tau(\TT)$ such that $(\tau(\TT))^\circ = \TT$ and $Q_{\tau(\TT)} = Q_{\TT}$ under the 
identification between $A(\TT)$ and $A_{\bowtie}(\tau(\TT))$.

It was proved in \cite{FST,FT} that cluster variables of $\AA(S,M)$ are indexed by $A(S,M)$ if $(S, M)$ is a once-punctured closed surface, and 
by  $A_{\bowtie}(S,M)$ otherwise. In the latter case the cluster variables will be denoted by $\{ x_\al \mid \al \in A_{\bowtie}(S,M)\}$ 
extending the notation from \S 2.2. Furthermore, \cite{FST,FT} proved that all seeds of $\AA(S,M)$ are 
$( \{ x_\al \mid \al \in A(\TT) \}, Q_\TT )$ for ordinary triangulations $\TT$ in the former case and 
$( \{ x_\al \mid \al \in A_{\bowtie}(\TT) \}, Q_\TT )$ for tagged triangulations $\TT$ in the latter case.
\sectionnew{Background on isomorphisms of cluster algebras}
\lb{backgr-auto}
In this section we gather background material on cluster isomorphisms and cluster automorphism groups.
\subsection{Cluster isomorphisms and cluster automorphism groups.}
\label{2b.1}
For a quiver $Q$, we will denote by $Q^{\opp}$ the opposite quiver.

Let $\AA_1$ and $\AA_2$ be cluster algebras of skewsymmetric geometric type without frozen variables. 
A {\em{strong cluster isomorphism}} \cite{FZ2} between $\AA_1$ and $\AA_2$ is an algebra isomorphism that 
maps a seed of $\AA_1$ to a seed of $\AA_2$. Such a map transports any seed of $\AA_1$ to a seed of $\AA_2$.
More formally:

\bde{sisom} \cite{FZ2} A strong cluster isomorphism between the cluster algebras $\AA_1$ and $\AA_2$ 
is an algebra isomorphism $\phi \colon \AA_1 \stackrel{\cong}\lra \AA_2$ such that for one
seed $({\bf{x}}, Q({\bf{x}}))$ of $\AA_1$, and thus for any seed of $\AA_1$, $(\phi({\bf{x}}), \phi(Q({\bf{x}})))$ is a seed of $\AA_2$.
Here $\phi(Q({\bf{x}}))$ denotes the quiver with vertex set $\phi({\bf{x}})$ which is isomorphic to $Q$ under the bijection of the vertex sets
${\bf{x}} \cong \phi({\bf{x}})$.
\ede

\bde{isom} \cite{ASS} A {\em{cluster isomorphism}} between $\AA_1$ and $\AA_2$ is 
an algebra isomorphism $\phi \colon \AA_1 \stackrel{\cong}\lra \AA_2$ such that for one
cluster ${\bf{x}}$ of $\AA_1$, and thus for any cluster of $\AA_1$: 
\begin{itemize}
\item $\phi({\bf{x}})$ is a cluster of $\AA_2$ and
\item for all $x \in {\bf{x}}$, $\mu_{x, {\bf{x}}}(x) = \mu_{\phi(x), \phi({\bf{x}})} ( \phi(x) )$
\end{itemize}
where $\mu_{x, {\bf{x}}}(x)$ denotes the mutation of the cluster variable $x \in {\bf{x}}$ in the direction of $x$.
\ede
\ble{isom2} \cite{ASS} Assume that $\AA_1$ and $\AA_2$ are two cluster algebras such that the exchange 
quiver of one seed of $\AA_1$ or $\AA_2$ (and thus of any seed of that cluster algebra) is connected.

An algebra isomorphism $\phi \colon \AA_1 \stackrel{\cong}\lra \AA_2$ is a cluster isomorphism
that has the property that for one seed $({\bf{x}}, Q({\bf{x}}))$ of $\AA_1$, and thus for any seed of $\AA_1$, 
either $(\phi({\bf{x}}), \phi(Q({\bf{x}})))$ or $(\phi({\bf{x}}), \phi(Q({\bf{x}}))^{\opp})$
is a seed of $\AA_2$. For such a map the two possibilities occur uniformly for all seeds of $\AA_1$.
\ele
Denote by $\Iso^+(\AA_1, \AA_2)$ and $\Iso(\AA_1, \AA_2)$ the sets of all strong cluster isomorphisms and 
all cluster isomorphisms, respectively. Set $\Aut^+(\AA_1):= \Iso^+(\AA_1, \AA_1)$ and $\Aut(\AA_1) :=\Iso(\AA_1, \AA_1)$
for the corresponding (strong) cluster automorphism groups. 

\bre{upper-aut}
The {\em{upper cluster algebra}} associated to the quiver $Q$ is defined \cite{BFZ} by 
\[
\UU({\bf{x}}, Q) = \bigcap_{\mbox{all clusters} \; \; \wt{\bf{x}} } \Zset [ x^{\pm 1}, x \in \wt{\bf{x}}].
\]
By the Laurent phenomenon \cite{FZ-l}, it contains the cluster algebra $\UU({\bf{x}}, Q)$. One can define the group of cluster automorphisms $\Aut  \UU({\bf{x}}, Q)$ as the set of all automorphisms of $\UU({\bf{x}}, Q)$ that satisfy the two conditions in \deref{isom}. It is easy to prove that there is a canonical isomorphism defined by restricting automorphisms of $\UU({\bf{x}}, Q)$ to $\AA({\bf{x}}, Q)$:
\[
\Aut \AA({\bf{x}}, Q) \cong \Aut \UU({\bf{x}}, Q) ,
\]
\ere
\subsection{Induced isomorphisms between surface cluster algebras.}
\label{2b.2}
A homeomorphism $g \colon (S_1, M_1) \stackrel{\cong}\lra (S_2, M_2)$ between two bordered 
(oriented) surfaces with marked points is a homeomorphism $g \colon S_1 \stackrel{\cong}\lra S_2$ such that $g(M_1) = M_2$.
Each homeomorphism $g \colon (S_1, M_1) \stackrel{\cong}\lra (S_2, M_2)$ gives rise to a 
cluster isomorphism 
\[
\psi_g \in \Iso ( \AA(S_1,  M_1), \AA(S_2, M_2))  \quad \mbox{defined by} \quad
\psi_g (x_{\al}) := x_{g(\al)} 
\]
for all $\al \in A_{\bowtie}(S_1, M_1)$. Furthermore, $\psi_g \in \Iso^+ ( \AA(S_1,  M_1), \AA(S_2, M_2))$ if and only if 
$g$ is an orientation-preserving homeomorphism.

Assume that $(S_2, M_2)$ is not a once-punctured closed surface. Let $R$ be a subset of the set of punctures of $(S_2, M_2)$. For a tagged arc 
$\al \in \AA_{\bowtie}(S_2, M_2)$, denote by $\al^R$ the tagged arc obtained from it by changing the taggings of $\al$ at those 
of its endpoints that belong to $R$. For a tagged triangulation $\TT$ of $(S_2,M_2)$, $\TT^R$ is also a tagged triangulation of $(S_2,M_2)$
having the same exchange quiver.
Fix a tagged triangulation $\TT$ of $(S_2, M_2)$ and define \cite[Lemma 4.9]{ASS}
\[
\psi_R \in \Aut^+ \AA(S_2, M_2)  \quad \mbox{by} \quad \psi_R(x_\al) = x_{\al^R}
\]
for $\al \in A_{\bowtie}(\TT)$ and extended to $\AA(S_2, M_2) \subset \Qset(x_\al, \al \in A_{\bowtie}(\TT) )$ by the algebra homomorphism property.
The map $\psi_R$ is independent of the choice of $\TT$ used to define it. This can be seen directly from the properties of mutations of tagged triangulations
\cite[\S 9.1-9.3]{FST}. Another easy justification is to argue that for every arc $x$ of a triangulation $\TT$, $\TT^R$ and $(\mu_x(\TT))^R$ differ by exactly one 
arc, and therefore must be a one-step mutation from each other with respect to this arc since $\psi_R$ is a strong cluster automorphism.

For $g$ and $R$ as above, set
\begin{equation}
\label{ind-clisom}
\psi_{g, R} := \psi_R \psi_g \in \Iso ( \AA(S_1,  M_1), \AA(S_2, M_2)).
\end{equation}
\subsection{Induced automorphisms of surface cluster algebras}
\label{2b.3}
We restrict the discussion from the previous subsection to the case when $S_2=S_1$ and $M_2 = M_1$.
For an oriented surface $S$ with boundary $\partial S$ denote by $\Homeo^+(S, \partial S)$ the group 
of orientation-preserving homeomorphisms of $S$ that fix $\partial S$. Let $\Homeo^+_0(S, \partial S)$ 
be its subgroup of those homeomorphisms that are isotopic to the identity. 
The {\em{mapping class group}} of $S$ is the factor group 
\[
\Mod S := \Homeo^+(S, \partial S) / \Homeo^+_0(S, \partial S).
\]

Analogously, for a bordered (oriented) surface with marked points $(S,M)$ denote by $\Homeo^+(S, M)$ 
the group of orientation-preserving homeomorphisms of $S$ that take $M$ to itself. (Such maps are not required to fix $M$ or $\partial S$.)
Let $\Homeo^+_0(S, M)$  be its subgroup of those homeomorphisms that are isotopic to the identity through an isotopy that pointwisely fixes $M$ at all times. 
The {\em{mapping class group}} 
of the pair $(S,M)$ is defined by 
\[
\MCG(S,M):= \Homeo^+(S,M) / \Homeo^+_0(S,M).
\]
Define \cite{ASS} the {\em{tagged mapping class group}} of $(S,M)$ to be the semidirect product
\begin{equation}
\label{tagMCG}
\MCG_{\bowtie} (S,M) := \MCG(S,M) \ltimes \Zset_2^P,
\end{equation}
where $P$ is the set of punctures of $(S,M)$. Identify $\Zset_2^P$ with the set of subsets of $P$  
by sending $R \subseteq P$ to the element of $\Zset_2^P$ with 1's in the positions of $R$ and 0's elsewhere. 
The semidirect product in \eqref{tagMCG} is defined with respect to the 
action of $\MCG(S,M)$ on $\Zset_2^P$ given by $g \cdot R := g(R)$. 
Assem, Schiffler, and Shramchenko, proved \cite{ASS} that the map $g \in \MCG(S,M)\mt \psi_g$
defines an embedding
\begin{equation}
\label{ASS-emb1}
\MCG(S,M) \hra \Aut^+ \AA(S,M) 
\end{equation}
for all $(S,M)$ and that the 
map $(g, R) \in \MCG_{\bowtie}(S,M) \mt \psi_{g,R}$ defines an embedding
\begin{equation}
\label{ASS-emb2}
\MCG_{\bowtie}(S,M) \hra \Aut^+ \AA(S,M) 
\end{equation}
if $(S,M)$ is not a once-punctured closed surface. 

We will also need an unsigned version of the mapping class group $\MCG(S,M)$. Denote by
$ \Homeo(S,M)$ the group of homeomorphisms of $S$ that take $M$ to itself. Define 
\begin{align*}
\MCG^\pm(S,M) &:= \Homeo(S,M) / \Homeo^+_0(S,M) \quad \mbox{and} 
\\
\MCG_{\bowtie}^\pm(S,M) &:= \MCG^\pm(S,M) \ltimes \Zset_2^P.
\end{align*} 
The above two embeddings imply that in the settings of \eqref{ASS-emb1} and \eqref{ASS-emb2}, respectively, 
we have the embeddings
\begin{align}
&\MCG^\pm(S,M) \hra \Aut \AA(S,M), \quad
g \in \MCG^\pm(S,M) \mt \psi_g,
\label{ASS-emb1u}
\\
&\MCG_{\bowtie}^\pm(S,M) \hra \Aut \AA(S,M), \quad
(g, R) \in \MCG_{\bowtie}(S,M) \mt \psi_{g,R}.
\label{ASS-emb2u}
\end{align}
\sectionnew{Maximal triangulations}
\lb{match-fa}
In this section we define and study maximal triangulations.
Two equivalent definitions of this notion are given in \deref{max} and \coref{max-triang}; \S \ref{3.1} contains 
auxiliary results needed to establish the equivalence. Two matching properties for maximal triangulations with isomorphic 
exchange quivers are proved in Propositions \ref{pmatch1} and \ref{pmatch2}. 
\subsection{Double-glued triangles and number of edges of $Q_\TT$} 
\label{3.1}
Let $(S,M)$ be a bordered surface with marked points and $\TT$ be a triangulation of $(S,M)$. Because of the exclusion of the sphere with 3 punctures 
in \deref{bor}, two triangles of $\TT$ cannot have 3 common sides unless $(S,M)$ is the once-punctured torus. Next, we consider the situation when a pair 
of triangles has exactly two common arcs; we call those {\em{double-glued}}.

\bde{double-glue} (a) A pair of triangles of a triangulation $\TT$ will be called {\em{positively double-glued}} if they have exactly two common arcs and those arcs 
appear in the same order when the boundaries of the triangles are traced clockwise. 

(b) A pair of triangles of $\TT$ will be called {\em{negatively double-glued}} if they have exactly 2 common arcs, but those arcs 
appear in the opposite order in the clockwise tracing of boundaries of the triangles. 

(c) Denote by $d_{neg}(\TT)$ the number of pairs of negatively double-glued triangles of $\TT$. 
(A triangle of $\TT$ cannot appear in two different pairs of double-glued triangles because each arc is on the boundary of two triangles.)

(d) A pair of triangles of $\TT$ will be called {\em{single-glued}} if they share one common arc.  
\ede
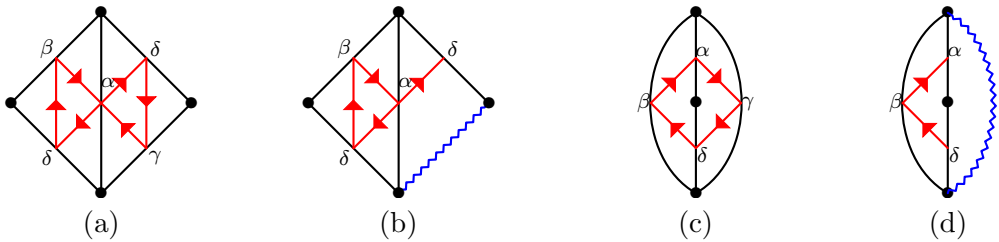
\begin{figure}[ht] 
\begin{center} 
\begin{tabular}{cccc}
\begin{tikzpicture}[line cap=round,line join=round,>=triangle 45,x=1.0cm,y=1.0cm,scale=.6]
\begin{scope}[thick, every node/.style={sloped,allow upside down}]
\clip(3,3.9) rectangle (9,9);
\draw (6,8)--(4,6);
\draw (6,8)--(8,6);
\draw (6,8)--(6,4);
\draw (4,6)--(6,4);
\draw (6,4)--(8,6);

\draw(4,6) node{$\bullet$};
\draw(6,8) node{$\bullet$};
\draw(6,4) node{$\bullet$};
\draw(8,6) node{$\bullet$};

\draw[red] (5,7) -- node{\midarrow} (6,6);
\draw[red] (6,6) -- node{\midarrow} (5,5);
\draw[red] (5,5) -- node{\midarrow} (5,7);

\draw[red] (6,6) -- node{\midarrow} (7,7);
\draw[red] (7,7) -- node{\midarrow} (7,5);
\draw[red] (7,5) -- node{\midarrow} (6,6);

\begin{scriptsize}
\draw(6.15,6.45) node {$\al$};
\draw(7.2,4.8) node {$\ga$};
\draw(7.2,7.2) node {$\de$};
\draw(4.8,4.8) node {$\de$};
\draw(4.8,7.2) node {$\be$};
\end{scriptsize}
\end{scope}
\end{tikzpicture}
&
\begin{tikzpicture}[line cap=round,line join=round,>=triangle 45,x=1.0cm,y=1.0cm,scale=.6]
\begin{scope}[thick, every node/.style={sloped,allow upside down}]
\clip(3,3.9) rectangle (9,9);
\draw (6,8)--(4,6);
\draw (6,8)--(8,6);
\draw (6,8)--(6,4);
\draw (4,6)--(6,4);
\draw[color=blue, line join=round, decorate, decoration={zigzag,
    segment length=4,
    amplitude=.9,post=lineto,
    post length=2pt
}](6,4)--(8,6);

\draw(4,6) node{$\bullet$};
\draw(6,8) node{$\bullet$};
\draw(6,4) node{$\bullet$};
\draw(8,6) node{$\bullet$};

\draw[red] (5,7) -- node{\midarrow} (6,6);
\draw[red] (6,6) -- node{\midarrow} (5,5);
\draw[red] (5,5) -- node{\midarrow} (5,7);

\draw[red] (6,6) -- node{\midarrow} (7,7);

\begin{scriptsize}
\draw(6.15,6.45) node {$\al$};
\draw(7.2,7.2) node {$\de$};
\draw(4.8,4.8) node {$\de$};
\draw(4.8,7.2) node {$\be$};
\end{scriptsize}
\end{scope}
\end{tikzpicture}
&
\begin{tikzpicture}[line cap=round,line join=round,>=triangle 45,x=1.0cm,y=1.0cm,scale=.6]
\begin{scope}[thick, every node/.style={sloped,allow upside down}]
\clip(1,1.9) rectangle (7,8);

\draw(4,2) node{$\bullet$};
\draw(4,4) node{$\bullet$};
\draw(4,6) node{$\bullet$};

\begin{scriptsize}
\draw(4.15,5.15) node {$\al$};
\draw(4.15,2.85) node {$\de$};
\draw(5.15,4) node {$\ga$};
\draw(2.85,4) node {$\be$};
\end{scriptsize}

\draw (4,2) -- (4,4);
\draw (4,4) -- (4,6);

\draw (4,2) edge[bend left=60] (4,6);
\draw (4,2) edge[bend right=60] (4,6);

\draw[red] (4,5) -- node{\midarrow} (5,4);
\draw[red] (5,4) -- node{\midarrow} (4,3);

\draw[red] (4,3) -- node{\midarrow} (3,4);
\draw[red] (3,4) -- node{\midarrow} (4,5);

\end{scope}
\end{tikzpicture}
&
\begin{tikzpicture}[line cap=round,line join=round,>=triangle 45,x=1.0cm,y=1.0cm,scale=.6]
\begin{scope}[thick, every node/.style={sloped,allow upside down}]
\clip(2,1.9) rectangle (6,8);

\draw(4,2) node{$\bullet$};
\draw(4,4) node{$\bullet$};
\draw(4,6) node{$\bullet$};

\begin{scriptsize}
\draw(4.15,5.15) node {$\al$};
\draw(4.15,2.85) node {$\de$};
\draw(2.85,4) node {$\be$};
\end{scriptsize}

\draw (4,2) -- (4,4);
\draw (4,4) -- (4,6);

\draw (4,2) edge[bend left=60] (4,6);
\draw (4,2) edge[color=blue, line join=round,bend right=60, decorate, decoration={zigzag,
    segment length=4,
    amplitude=.9,post=lineto,
    post length=2pt
}]  (4,6);


\draw[red] (4,3) -- node{\midarrow} (3,4);
\draw[red] (3,4) -- node{\midarrow} (4,5);

\end{scope}
\end{tikzpicture}\\
(a)&(b)&(c)&(d)\\

\end{tabular}
\end{center}
\caption{Double-glued triangles
\label{doub-glue}}
\end{figure}
\noindent
Recall \deref{triangles}. There are two types of pairs of double-glued triangles: 
\begin{enumerate}
\item Two (positively/negatively) double-glued faces, as on \figref{doub-glue} (a) and (c), respectively.
\item A (positively/negatively) double-glued pair of a face and a wedge, as on \figref{doub-glue} (b) and (d), 
respectively.
\end{enumerate}
{\bf{Convention for cutouts.}} As on \figref{extriang}, solid black lines indicate arcs and wavy blue lines indicate boundary arcs. Here and below the red arrows 
denote the edges of the exchange quiver $Q_\TT$ before we remove two cycles from the quiver. 
The cutouts on \figref{doub-glue} (a) and (c) have to be glued along the arc $\de$. There is a unique way in which 
this can be done so that the resulting subsurface is oriented. For any of the cutouts in the paper there will always be a unique way to glue them along the repeated arcs so that the resulting subsurface is oriented.
\medskip
\\
The difference between a positively and negatively double-glued pair of triangles is that in the latter case a 2-cycle is removed 
in the construction of the exchange quiver $Q_\TT$, recall \S \ref{2.2}. In the former case we have no removal of cycles, instead there 
is a double arrow between the vertices of $Q_\TT$ corresponding to the 2 common arcs.

\begin{figure}[ht] 
\begin{center} 
\begin{tabular}{ccc}
\begin{tikzpicture}[scale=.7,baseline=0pt]
\draw (0,0) .. controls (1,2) and (-1,2) .. (0,0);
\draw (0,0) .. controls (1,1) and (1,2) .. (0,3);
\draw (0,0) .. controls (-1,1) and (-1,2) .. (0,3);
\draw (0,0) node {$\bullet$}; \draw (0,3) node {$\bullet$};
\draw (0,0) -- (0,1);
\draw (0,1) node {$\bullet$}; 
\draw (2,0) node {};
\draw (-2,0) node {};
\end{tikzpicture} 
&
\begin{tikzpicture}[scale=.7,baseline=0pt]
\draw (0,0) .. controls (-2.7,3) and (0,3) .. (0,0);
\draw (0,0) .. controls (0,3) and (2.7,3) .. (0,0);
\draw (0,1.5) circle (1.5cm);
\draw (0,0) node {$\bullet$}; 
\draw (0,0) -- (.8,1.6);
\draw (0,0) -- (-.8,1.6);
\draw (.8,1.6) node {$\bullet$}; \draw (-.8,1.6) node {$\bullet$}; 
\end{tikzpicture} 
&
\begin{tikzpicture}[scale=.7,baseline=0pt]
\draw (0,0) .. controls (3,2) and (3,-2) .. (0,0);
\draw (0,0) .. controls (2,3) and (-2,3) .. (0,0);
\draw (0,0) .. controls (-3,2) and (-3,-2) .. (0,0);
\draw (0,0) node {$\bullet$}; 
\draw (0,0) -- (1,0);
\draw (0,0) -- (-1,0);
\draw (0,0) -- (0,1);
\draw (1,0) node {$\bullet$}; \draw (-1,0) node {$\bullet$}; \draw (0,1) node {$\bullet$}; 
\end{tikzpicture} \\
(a) & (b) & (c)\\

\end{tabular}
\end{center}
\caption{Triangles whose boundaries contain loops
\label{self-fold}}
\end{figure}
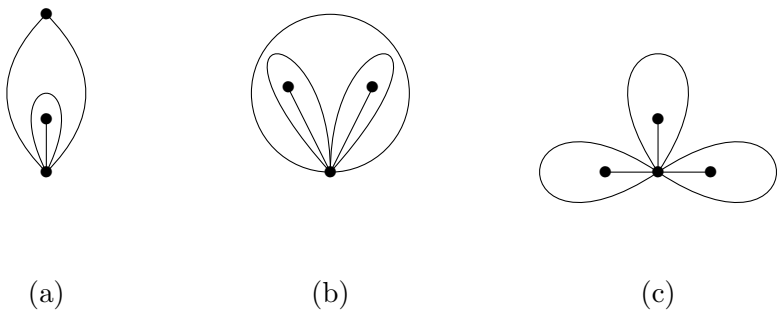

Consider the loop of a self-folded triangle $\De$ of $\TT$. It is in the boundary of exactly one more triangle $\De'$ which is not self-folded because the 
sphere with 3 punctures is excluded in \deref{bor}.
The boundary of $\De'$ can contain 1, 2 or 3 loops enclosing self-folded triangles as on Fig. \ref{self-fold} (a), (b), and (c), respectively. The third case is only possible for the sphere with 4 punctures.  
In general, the outer arcs in \figref{self-fold} (a) and (b) can be either ordinary or boundary arcs. 

Denote by $s_f(\TT)$ the number of different faces of a triangulation $\TT$ whose boundary contains at least one loop of a self-folded 
triangle. Denote by $s_w(\TT)$ the number of wedges of $\TT$ whose boundary contains a loop.

\bpr{edges} Let $\TT$ be a triangulation of a bordered surface with marked points $(S,M)$ which is different from the three
triangulations in \figref{self-fold} of the once-punctured digon, twice-puncture monogon, and 4-punctured sphere, respectively. Then the 
number of edges $e(Q_\TT)$ of $\TT$ is given by 
\[
e(Q_\TT) = 3 f(\TT) + w(\TT) - 2 d_{neg}(\TT) - s_f(\TT) - 2 s_w(\TT),
\] 
recall \S \ref{2.2} and \deref{double-glue} (c) for the definition of the numbers in the rhs.
\epr
In the language of Fomin, Shapiro and Thurston \cite{FST}, $s_f(\TT)$ counts puzzle pieces in $\TT$ of the type in \figref{self-fold} (a), (b) that 
do not contain a boundary arc; $s_w(\TT)$ counts the puzzle pieces in $\TT$ of the type in \figref{self-fold} (a), (b) that 
have one boundary arc. Because of the exclusions in the proposition, the outer triangles in \figref{self-fold} (a) and (b) cannot be 
a cap and a wedge, respectively. In terms of block decompositions of exchange quivers \cite{FST} (see \S \ref{4.2} for details on this construction) 
$s_f(\TT)$ counts the number of blocks of type IV and V in $Q_{\TT}$; $s_w(\TT)$ counts the number of blocks of type III
in $Q_{\TT}$.
\medskip
\\
\noindent
{\em{Proof of \prref{edges}.}} The proof can be phrased in terms of counting number of edges of block decompositions of exchange quivers
(see \S \ref{4.2}) or direct counting of the number of edges of $Q_\TT$ using its definition. We follow the latter:
In the first step of the computation of the exchange quiver $Q_\TT$ (see \S \ref{2.2}), if $\TT$ has no self-folded triangles, then 
each face of $\TT$ contributes 3 edges and each wedge of $\TT$ contributes 1 edge. 

If the puzzle pieces on Fig. \ref{self-fold} (a) and (b) contain only faces, then they contribute 
$8= 3 \cdot 3-1$ and $5= 2 \cdot 3 -1$ edges to the first step of the computation of $Q_\TT$, respectively.  
If the outer triangle on Fig. \ref{self-fold} (a) is a wedge, then the puzzle piece contributes 
$4= 2 \cdot 3 + 1 - 3$ edges to the first step of the computation of $Q_\TT$. Thus, in the first 
step of the computation of $Q_\TT$ we have 
\[
e(Q_\TT) = 3 f(\TT) + w(\TT) - s_f(\TT) - 2 s_w(\TT),
\]
faces. In the second part of the computation of $Q_\TT$, each pair $\De_1, \De_2$ of negatively double-glued triangles leads to the removal of exactly 
one 2-cycle between the 2 vertices of $Q_\TT$ which correspond to the 2 common arcs of $\De_1, \De_2$. Here we use the fact that none of the common 
arcs is a loop or a radius of a self-folded triangle because each arc is in the boundary of exactly 2 triangles. 
\qed
\subsection{Maximal triangulations}
\label{3.2}
\bde{max} We will say that a triangulation $\TT$ of $(S,M)$ is {\em{maximal}} if it has no self-folded triangles 
or pairs of negatively double-glued triangles, and the number of its caps is the maximal possible: 
$c(\TT) = c(S,M)$, recall \eqref{cSM} and \eqref{caps}. 
\ede

{\em{Any surface that is different from the the once-punctured 2-, 3- and 4-gons, 
and the twice-punctured monogon has a maximal triangulation}}; a stronger statement is proved in 
\thref{exist}. The once-punctured digon and the twice-punctured monogon do not have maximal triangulations because each of their triangulations 
either has a negatively double-glued pair of triangles or a self-folded triangle. Any triangulation of the once-punctured triangle and square with maximal number of caps will also contain a negatively double-glued pair of triangles or a self-folded triangle.

The next result shows that maximal triangulations are precisely those triangulations whose quivers have maximum number 
of arrows among all quivers in the mutation class, except for the 4-punctured sphere and the four surfaces that do not admit maximal 
triangulations.

\bco{max-triang} Assume that $\TT$ is a triangulation of a bordered surface with marked points $(S,M)$ 
which is different from the triangulation in \figref{self-fold} (c) of the 4-punctured sphere.
Then 
\[
e(Q_\TT) \leq 2 n(S,M) - h(S,M) + c(S,M)
\]
and equality holds if and only if $\TT$ is maximal.
\eco
For the four surfaces without maximal triangulations, the last part of the corollary means that the inequality is strict 
for all triangulations.
\begin{proof}
\prref{edges}, eqs. \eqref{fwc1}--\eqref{fwc2} and the inequality \eqref{caps} 
imply
\[
e(Q_\TT) \leq 3 f(\TT) + w(\TT) = 2n -h + c(\TT) \leq  2n - h + c(S,M).
\]
The two inequalities turn into equalities precisely when $\TT$ is maximal.
\end{proof}

If $S$ is a closed surface, the condition that a triangulation $\TT$ has no self-folded triangles or pairs of negatively double-glued 
triangles is equivalent to saying that each puncture is adjacent to at least 3 arcs. The latter condition was used in \cite[\S 2.1]{L}, 
however, the results in \cite[\S 2.1]{L} do not extend to the case of nonclosed surfaces.
\bre{4sph} The quiver of the triangulation in \figref{self-fold} (c) of the 4-punctured sphere has 12 edges. Combining this 
with \coref{max-triang}, we obtain the following:

{\em{The quiver of each triangulation of the 4-punctured sphere has $\leq 12$ edges and equality is obtained if and only if the triangulation is maximal or is 
the one in \figref{self-fold}}} ({\em{c}}) (with the various different labelings of the vertices).
\ere
\subsection{Matching maximal triangulations under quiver isomorphsims}
\label{3.3}
\bpr{match1} Let $(S_1,M_1)$ and $(S_2,M_2)$ be bordered surfaces with marked points 
which are different from the 4-punctured sphere, the twice-punctured monogon, and the once-punctured 2-, 3- and 4-gons.
Assume that $\TT_1$ and $\TT_2$ are tagged triangulations of $(S_1,M_1)$ and $(S_2,M_2)$, and that
there is an isomorphism
\[
\phi \colon Q_{\TT_1} \stackrel{\cong}{\lra} Q_{\TT_2}.
\]
If $(\TT_1)^\circ$ is maximal, then $(\TT_2)^\circ$ is also maximal.
\epr
\begin{proof} Since $Q_{\TT_i} \cong Q_{(\TT_i)^\circ}$, we can assume that $\TT_1$ and $\TT_2$ are ordinary triangulations.
By \coref{max-triang}, 
\begin{multline}
\label{ineq1}
2 n(S_1,M_1) - h(S_1,M_1) + c(S_1,M_1)
= e(Q_{\TT_1}) = \\
e(Q_{\TT_2}) \leq 2 n(S_2,M_2) - h(S_2,M_2) + c(S_2,M_2).
\end{multline}
Because $(S_2, M_2)$ is different from the surfaces listed above it has a maximal triangulation. Therefore there exists a sequence of mutations $s$ such that $\mu_s (\TT_2)$ is a maximal triangulation.
At the same time we have an induced isomorphism $ \mu_{\phi(s)} (Q_{\TT_1}) \stackrel{\cong}{\lra} \mu_s( Q_{\TT_2})$.
Applying again \coref{max-triang}, gives
\begin{multline}
\label{ineq2}
2 n(S_1,M_1) - h(S_1,M_1) + c(S_1,M_1)
\geq e( \mu_{\phi(s)} (Q_{\TT_1})) = \\ e(\mu_s(Q_{\TT_2})) =  2 n(S_2,M_2) - h(S_2,M_2) + c(S_2,M_2).
\end{multline}
This is only possible if both \eqref{ineq1} and \eqref{ineq2} are equalities, in which case \coref{max-triang} implies that 
$\TT_2$ has no self-folded triangles or pairs of negatively double-glued triangles.
\end{proof}
\subsection{Matching pairs of glued faces}
\label{3.4}
\bpr{match2} Let $(S_1,M_1)$ and $(S_2,M_2)$ be two bordered surfaces with marked points 
which are different from
\begin{enumerate}
\item the 4-punctured sphere,
\item the twice-punctured monogon and digon,  and
\item
the once-punctured 2-, 3- and 4-gons.
\end{enumerate}
Assume that $\TT_1$ and $\TT_2$ are triangulations of $(S_1,M_1)$ and $(S_2,M_2)$
such that $\TT_1$ is maximal and that
\[
\phi \colon Q_{\TT_1} \stackrel{\cong}{\lra} Q_{\TT_2} 
\]
is an isomorphism. 

(i) Let $\al_i, \be_i, \ga_i, \de_i, \zeta_i$ be distinct arcs of $\TT_i$ such that $\phi([\al_1]) = [\al_2]$, $\phi([\be_1]) = [\be_2]$, etc. If, 
$\{ \al_1, \be_1, \ga_1 \}$ and $\{\al_1, \de_1, \zeta_1 \}$ are the clockwise boundaries of two faces of $\TT_1$, then 
$\{ \al_2, \be_2, \ga_2 \}$ and $\{\al_2, \de_2, \zeta_2 \}$ are the clockwise boundaries of two faces of $\TT_2$.

(ii) Let $\al_i, \be_i, \ga_i, \de_i$ be distinct arcs of $\TT_i$ such that $\phi([\al_1]) = [\al_2]$, $\phi([\be_1]) = [\be_2]$, etc. If, 
$\{ \al_1, \de_1, \be_1 \}$ and $\{\al_1, \de_1, \ga_1 \}$ are the clockwise boundaries of two faces of $\TT_1$, then 
$\{ \al_2, \de_2, \be_2 \}$ and $\{\al_2, \de_2, \ga_2 \}$ are the clockwise boundaries of two faces of $\TT_2$.
\epr
Informally, the first part of the theorem states that, with the exclusion of the surfaces (1)-(3), every isomorphism of exchange graphs 
$\phi \colon Q_{\TT_1} \stackrel{\cong}{\lra} Q_{\TT_2} $, maps the vertices corresponding to the arcs of 
2 single-glued faces of $\TT_1$ to the vertices corresponding to the arcs of 
2 single-glued faces of $\TT_2$. The second part of the theorem states the same thing for positively double-glued pairs 
of faces.
\medskip
\\
{\em{Proof of \prref{match2}.}} In each of parts (a) and (b), $Q_{\TT_2}$ needs to have a 3-cycle, 
so $(S_2, M_2)$ should be different from the once-punctured digon.
\prref{match1} implies that $\TT_2$ is a maximal triangulation. 

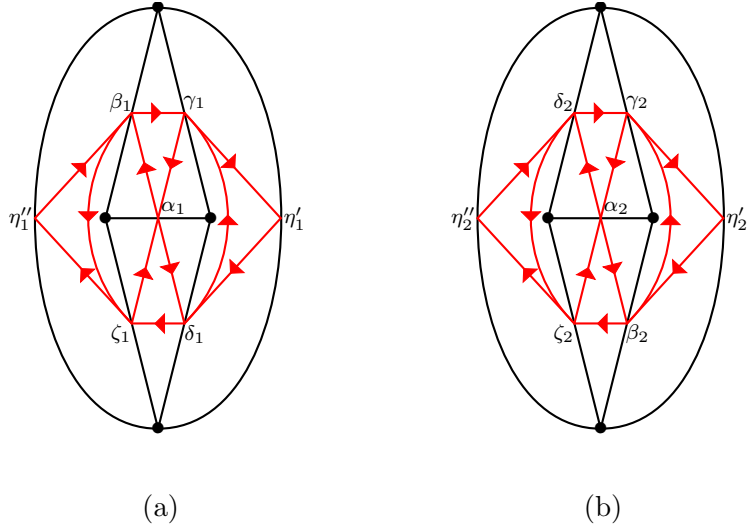
\begin{figure}[ht]
\begin{center}
\begin{tabular}{cc}
\begin{tikzpicture}[line cap=round,line join=round,>=triangle 45,x=1.0cm,y=1.0cm,scale=.7]
\begin{scope}[thick, every node/.style={sloped,allow upside down}]
\clip(4.1,3) rectangle (12,12.1);
\draw (8,4) -- (9,8);
\draw (9,8) -- (8,12);
\draw (8,4) -- (7,8);
\draw (7,8) -- (8,12);
\draw (7,8) -- (9,8);
\draw (8,4) edge[bend left=90] (8,12);
\draw (8,4) edge[bend right=90] (8,12);

\draw (8,4) node {$\bullet$};
\draw (9,8) node {$\bullet$};
\draw (7,8) node {$\bullet$};
\draw (8,12) node {$\bullet$};

\draw[red] (10.33,8) -- node{\midarrow} (8.5,6);
\draw[red] (8.5,10) -- node{\midarrow} (10.33,8);
\draw[red] (8.5,6) edge[bend right=45] node{\midarrow} (8.5,10);

\draw[red] (7.5,6) -- node{\midarrow} (5.66,8);
\draw[red] (5.66,8) -- node{\midarrow} (7.5,10);
\draw[red] (7.5,10) edge[bend right=45] node{\midarrow} (7.5,6);

\draw[red] (7.5,10) -- node{\midarrow} (8.5,10);
\draw[red] (8.5,10) -- node{\midarrow} (8,8);
\draw[red] (8,8) -- node{\midarrow} (7.5,10);

\draw[red] (8,8) -- node{\midarrow} (8.5,6);
\draw[red] (8.5,6)-- node{\midarrow} (7.5,6);
\draw[red] (7.5,6)-- node{\midarrow} (8,8);
\begin{scriptsize}
\draw(10.6,8) node {$\eta'_1$};
\draw(5.4,8) node {$\eta''_1$};
\draw(8.3,8.2) node {$\al_1$};
\draw(8.7,10.2) node {$\ga_1$};
\draw(7.3,10.2) node {$\be_1$};
\draw(8.7,5.8) node {$\de_1$};
\draw(7.3,5.8)node {$\zeta_1$};
\end{scriptsize}
\end{scope}
\end{tikzpicture}
&
\begin{tikzpicture}[line cap=round,line join=round,>=triangle 45,x=1.0cm,y=1.0cm,scale=.7]
\begin{scope}[thick, every node/.style={sloped,allow upside down}]
\clip(4.1,3) rectangle (12,12.1);

\draw (8,4) -- (9,8);
\draw (9,8) -- (8,12);
\draw (8,4) -- (7,8);
\draw (7,8) -- (8,12);
\draw (7,8) -- (9,8);
\draw (8,4) edge[bend left=90] (8,12);
\draw (8,4) edge[bend right=90] (8,12);

\draw (8,4) node {$\bullet$};
\draw (9,8) node {$\bullet$};
\draw (7,8) node {$\bullet$};
\draw (8,12) node {$\bullet$};

\draw[red] (10.33,8) -- node{\midarrow} (8.5,6);
\draw[red] (8.5,10) -- node{\midarrow} (10.33,8);
\draw[red] (8.5,6) edge[bend right=45] node{\midarrow} (8.5,10);

\draw[red] (7.5,6) -- node{\midarrow} (5.66,8);
\draw[red] (5.66,8) -- node{\midarrow} (7.5,10);
\draw[red] (7.5,10) edge[bend right=45] node{\midarrow} (7.5,6);

\draw[red] (7.5,10) -- node{\midarrow} (8.5,10);
\draw[red] (8.5,10) -- node{\midarrow} (8,8);
\draw[red] (8,8) -- node{\midarrow} (7.5,10);

\draw[red] (8,8) -- node{\midarrow} (8.5,6);
\draw[red] (8.5,6)-- node{\midarrow} (7.5,6);
\draw[red] (7.5,6)-- node{\midarrow} (8,8);
\begin{scriptsize}
\draw(10.6,8) node {$\eta'_2$};
\draw(5.4,8) node {$\eta''_2$};
\draw(8.3,8.2) node {$\al_2$};
\draw(8.7,10.2) node {$\ga_2$};
\draw(7.3,10.2) node {$\de_2$};
\draw(8.7,5.8) node {$\be_2$};
\draw(7.3,5.8)node {$\zeta_2$};
\end{scriptsize}
\end{scope}
\end{tikzpicture}\\
(a)&(b)\\
\end{tabular}
\end{center}
\caption{Quiver isomorphisms for single-glued triangles
\label{big}}
\end{figure}

(i) The quiver $Q_{\TT_2}$ has the 4 edges $[\al_2] \to [\be_2], [\de_2]$ and $[\ga_2], [\zeta_2] \to [\al_2]$. 
Assume that the statement is not correct. Then the 2 triangles of $\TT_2$ whose boundaries contain $\al_2$ are faces 
and have clockwise boundaries $\{ \al_2, \de_2, \ga_2 \}$ and $\{\al_2, \be_2, \zeta_2 \}$, as on \figref{big} (b).
It follows from the existence of the isomorphism $\phi \colon Q_{\TT_1} \to Q_{\TT_2}$ that 
$Q_{\TT_1}$ has edges $[\be_1] \to [\zeta_1]$ and $[\de_1] \to [\ga_1]$, 
and $Q_{\TT_2}$ has edges $[\de_2] \to [\zeta_2]$ and $[\be_2] \to [\ga_2]$. 
Therefore, $\TT_1$ has faces whose clockwise boundaries are  $\{\ga_1, \de_1, \eta'_1\}$ and 
$\{\be_1, \zeta_1, \eta''_1\}$ as on \figref{big} (a). Similarly, 
 $\TT_2$ should have faces whose clockwise boundaries are
 $\{\be_2, \ga_2, \eta'_2\}$ and $\{\de_2, \zeta_2, \eta''_2\}$.

 At least one of $\eta'_1, \eta'_1, \eta''_2, \eta''_2$
 should not be a boundary arc. Say $\eta'_1$ is an arc; the other cases are analogous.
 The quiver $Q_{\TT_1}$ has the 2 edges $[\ga_1] \to [\al_1], [\eta'_1]$. The quiver
 $Q_{\TT_2}$ has the edge $[\ga_2] \to [\al_2]$ and the only other outgoing edge 
 from $[\ga_2]$ that it can have is $[\ga_2] \to [\eta'_2]$ because the 2 triangles whose boundaries
 contain $\ga_2$ are already identified. Thus, $\eta'_2$ is an arc and $\phi([\eta'_1]) = [\eta'_2]$.
 This argument can now be repeated. The quiver $Q_{\TT_2}$ has the 2 edges $[\al_2], [\eta'_2] \to [\be_2]$. The quiver
 $Q_{\TT_1}$ has the edge $[\al_1] \to [\be_1]$ and the only other incoming edge 
 to $[\be_1]$ that it can have is $[\eta''_1] \to [\be_1]$. Therefore, $\eta''_1$ is an arc and $\phi([\eta''_1]) = [\eta'_2]$.  

Finally,  the quiver $Q_{\TT_1}$ has the 2 edges $[\al_1], [\eta'_1] \to [\de_1]$. The quiver
$Q_{\TT_2}$ has the edge $[\al_2] \to [\de_2]$ and the only other incoming edge 
to $[\de_2]$ that it can have is $[\eta''_2] \to [\de_2]$. So, $\eta''_2$ should be an arc and $\phi([\eta'_1]) = [\eta''_2]$.  
 
Since, $\phi([\eta'_1]) = [\eta'_2]$, $\phi([\eta''_1]) = [\eta'_2]$ and $\phi([\eta'_1]) = [\eta''_2]$,
\[
\eta'_1 = \eta''_1 \quad \mbox{and} \quad \eta'_2 =\eta''_2. 
\]
This implies that $(S_1, M_1)$ and $(S_2, M_2)$ are both isomorphic to the 4-punctured sphere, which is a 
contradiction.

(ii) Consider the two triangles $\De', \De''$ of $\TT_2$ whose boundaries contain $\al_2$.
The quiver $Q_{\TT_2}$ must have the edges $[\be_2], [\ga_2] \to [\al_2]$ and the double edge $[\al_2] \Rightarrow [\de_2]$.
Thus, the clockwise boundaries of $\De'$ and $\De''$ are $\{ \al_2, \de_2, \be_2 \}$ and $\{\al_2, \de_2, \ga_2 \}$, 
which implies the second statement of the theorem. 
\qed
\sectionnew{Recovering topology from maximal triangulations}
\label{homeo} 
We will call a maximal triangulation $\TT$ of a bordered surface with marked points $(S,M)$ {\em{connected}} if 
\begin{equation}
\label{F(T)}
F(\TT):=\bigcup \{ \mbox{faces of $\TT$} \} \backslash M 
\end{equation}
is connected. This is equivalent to saying that the subgraph of $Q_\TT$ with vertices corresponding to the sides of all faces and 
those edges of $Q_\TT$ that come from faces, is connected. 

In \S \ref{4.1} we prove the main result in the paper, a generalization of Theorem A for tagged triangulations.
The theorem is equivalent to a uniqueness statement for the block decompositions of the 
exchange quivers of connected maximal triangulations with at least 2 faces; 
this is discussed in \S \ref{4.2}. In \S \ref{4.3} we classify the bordered surfaces with marked points $(S,M)$ that posses
connected maximal triangulations with at least 2 faces.
\subsection{Recovering topology from cluster structure for connected maximal triangulations with at least 2 faces}
\label{4.1}
For two subsets $X$ and $Y$ of a given set, denote by $X \ominus Y = (X \cup Y) \backslash (X \cap Y)$ their symmetric 
difference.

\bth{homeo}  Let $(S_1, M_1)$ and $(S_2, M_2)$ be two 
bordered surfaces with marked points which are different 
from the 4-punctured sphere, and the twice-punctured monogon and digon.

Assume that $\TT_1$ and $\TT_2$ are tagged triangulations of $(S_1, M_1)$ and $(S_2, M_2)$ such that $\TT^\circ_1$ is  
a connected maximal triangulation with at least 2 faces, and that
$\phi \colon Q_{\TT_1} \stackrel{\cong}\lra Q_{\TT_2}$ is an isomorphism of quivers. {\em{(}}The last condition is exactly the same as 
specifying a strong cluster isomorphism $\phi \in \Iso^+( \AA(S_1, M_1), \AA(S_2, M_2))$.{\em{)}} Then there exist an orientation-preserving 
homeomorphism $g \colon (S_1, M_1) \stackrel{\cong}\lra (S_2, M_2)$ such that
\[
\phi = \psi_{g, R_2 \ominus g(R_1)}
\]
where $R_i$ are subsets of those punctures of $(S_i, M_i)$ that have index $-1$ with respect to $\TT_i$. 
In the right hand side, the notation \eqref{ind-clisom} is used with respect to the triangulation $\TT_2^\circ$ of $(S_2,M_2)$.  
\eth
The surface $(S_2, M_2)$ should be different from the once-punctured 2-, 3- and 4-gons because the quivers of 
all triangulations of those surfaces have fewer than two 3-cycles. This implies that the surface $(S_2, M_2)$ should have a 
maximal triangulation, and by  \prref{match1}, $\TT^\circ_2$ should be such. 
Denote the induced isomorphism
\[
\phi^\circ \colon Q_{\TT^\circ_1} \stackrel{\cong}\lra Q_{\TT^\circ_2}.
\] 
By abuse of notation we will denote by the same letter the corresponding strong cluster 
isomorphism $\phi^\circ \in \Iso^+( \AA(S_1, M_1), \AA(S_2, M_2))$. 

The statement is straightforward to verify when $(S_1,M_1)$ is the once-punctured torus (using the fact that this is the 
only surface with a pair of triple-glued faces). For the remainder of this subsection, we will assume that 
$(S_1, M_1)$ and $(S_2, M_2)$ are not the once-punctured torus, so each $\TT_i$ does not have two triangles 
with exactly the same boundaries. 

\ble{aux} In the setting of \thref{homeo}, there exists an orientation preserving homeomorphism $g \colon (S_1, M_1) \stackrel{\cong}\lra (S_2, M_2)$ 
such that $\phi^\circ = \psi_g$.
\ele
\begin{proof} We will construct such a map $g$ in three stages:

(1) Construct a homeomorphism $g \colon F(\TT_1) \stackrel{\cong}\lra F(\TT_2)$ satisfying
\begin{equation}
\label{ma}
g(\al_1 \backslash M) = \al_2 \backslash M \quad \mbox{for all $\al_i \in A(\TT_i)$ such that $\phi([\al_1]) = [\al_2]$}
\end{equation} 
for those $\al_1, \al_2$ that belong to the boundaries of some faces of $\TT_1, \TT_2$, recall \eqref{F(T)}.

(2)  Construct an extension 
\begin{equation}
\label{g-inter}
g \colon \bigcup \{ \mbox{faces and wedges of $\TT_1$} \} \backslash M_1 \stackrel{\cong}\lra 
                                           \bigcup \{ \mbox{faces and wedges of $\TT_2$} \} \backslash M_2
\end{equation}
of the  homeomorphism from (1) which satisfies \eqref{ma} and $g(\partial S_1 \backslash M_1) = \partial S_2 \backslash M_2$
for the part of the boundary of $S_1$ that comes from wedges (and not caps).

(3) Construct an extension $g \colon S_1 \stackrel{\cong}\lra S_2$ of the  homeomorphism
from (2) which satisfies \eqref{ma}, $g(M_1) = M_2$ and $g(\partial S_1) = \partial S_2$.

Part (1): By interchanging the roles of $(S_1, M_1)$ and $(S_2, M_2)$ 
we can assume that $f(\TT_1) \geq f(\TT_2)$. Enumerate the faces of $\TT_1$: $\De^{(1)}_1, \ldots, \De^{(m)}_1$ in such a way 
that for $2 \leq k \leq m$, $\De^{(k)}_1$ shares a common arc with at least one of the faces $\De^{(1)}_1, \ldots, \De^{(k-1)}_1$. 
Inductively, we will identify distinct faces $\De^{(1)}_2, \ldots, \De^{(m)}_2$ of $\TT_2$ and construct homeomorphisms 
\begin{equation}
\label{g-step}
g \colon (\De^{(1)}_1 \cup \ldots \cup \De^{(k)}_1) \backslash M \stackrel{\cong}\lra (\De^{(1)}_2 \cup \ldots \cup \De^{(k)}_2) \backslash M
\end{equation}
with the property \eqref{ma} that extend the ones from the previous steps. In particular, this will imply that $f(\TT_1) = f(\TT_2)$.

The faces $\De^{(1)}_1$ and $\De^{(2)}_1$ have at least one common arc and are not negatively double-glued. \prref{match2} implies that 
there are two faces $\De^{(1)}_2$ and $\De^{(2)}_2$ of $\TT_2$ whose boundaries are matched under $\phi$ (as vertices of $Q_{\TT_1}$ and $Q_{\TT_2}$).
It follows from \prref{match2} that there exists 
a  homeomorphism $g \colon (\De^{(1)}_1 \cup \De^{(2)}_1) \backslash M \stackrel{\cong}\lra (\De^{(1)}_2 \cup \De^{(2)}_2) \backslash M$
satisfying \eqref{ma}. At the $k$-th step we have a face $\De^{(k)}_1$ sharing at least one common arc with a face $\De^{(j)}_1$ for some $j < k$ 
and such that $\De^{(k)}_1$ and $\De^{(j)}_1$ are not negatively double-glued. Applying \prref{match2}, associates a face $\De^{(k)}_2$ of $\TT_2$
-- this is the unique face whose clockwise boundary is $\{\al_2, \be_2, \ga_2\}$ and $\al_2, \be_2, \ga_2 \in A(\TT_2)$ are such that 
$\phi ([\al_1]) = [\al_2]$, $\phi([\be_1]) = [\be_2]$, $\phi([\ga_1]) = [\ga_2]$ and $\{\al_1, \be_1, \ga_1\}$ is the clockwise boundary of 
$\De^{(k)}_1$. There exists an orientation preserving homeomorphism $\De^{(k)}_1 \backslash M \stackrel{\cong}\lra \De^{(k)}_2 \backslash M$
satisfying \eqref{ma} that extends the restriction of the map $g$ from the previous step to the nonempty set of arcs
\[
(\partial \De^{(k)}_1 \backslash M_1) \cap (\De^{(1)}_1 \cup \ldots \cup \De^{(k-1)}_1).
\] 
This is the needed extension of $g$ for the $k$-th step. It completes the first part of the construction of $g$.

Part (2): We enumerate the wedges of $\TT_1$: $\De^{(m+1)}_1, \ldots, \De^{(l)}_1$
in such a way that for $m+1 \leq k \leq l$, $\De^{(k)}_1$ has a common arc with at least one of the triangles $\De^{(1)}_1, \ldots, \De^{(k-1)}_1$. 
Inductively, we will identify distinct wedges $\De^{(m+1)}_2, \ldots, \De^{(l)}_2$ of $\TT_2$ and construct homeomorphisms like in \eqref{g-step}
with the properties \eqref{ma} and $g (\partial S_1 \backslash M_1) = \partial S_2 \backslash M_2$ 
that extend the ones from the previous steps. Denote the clockwise boundary of $\De^{(k)}_1$ by 
$\{\al_1, \be_1, \eta_1\}$ where $\al_1, \be_1$ are distinct arcs and $\eta_1$ is a boundary arc. Denote by $\al_2$ and $\be_2$ the distinct 
arcs of $\TT_2$ such that $\phi([\al_1]) = [\al_2]$, $\phi([\be_1]) = [\be_2]$. Then $Q_{\TT_1}$ has an extra edge
$[\al_1] \to [\be_1]$ compared to its subquiver coming from $\De^{(1)}_1 \cup \ldots \cup \De^{(k-1)}_1$ (which by itself will 
contain one more edge $[\al_1] \to [\be_1]$ if $\De^{(k)}_1$ is positively double-glued to a face of $\TT_1$). Since the homeomorphism 
$g$ from the $(k-1)$-st step has the property \eqref{ma}, $Q_{\TT_2}$ has an extra edge
$[\al_2] \to [\be_2]$ compared to its subquiver coming from $\De^{(1)}_2 \cup \ldots \cup \De^{(k-1)}_2$.  
Because all faces of $\TT_2$ are already listed in $\De^{(1)}_1 \cup \ldots \cup \De^{(k-1)}_1$, 
$\TT_2$ has a wedge with clockwise boundary $\{\al_2, \be_2, \eta_2\}$ for some $\eta_2 \in BA(S_2,M_2)$.

At lest one of the arcs $\al_1$ and $\be_1$ belongs to $\De^{(1)}_1 \cup \ldots \cup \De^{(k-1)}_1$, and, by the 
properties of $g$ from the previous step, the exact the corresponding ones among the arcs $\al_2$ and $\be_2$ 
belong to $\De^{(1)}_2 \cup \ldots \cup \De^{(k-1)}_2$. As in part one, from the previous step of the construction of
$g$ we have a homeomorphism 
\[
g \colon (\partial \De^{(k)}_1 \backslash M_1) \cap (\De^{(1)}_1 \cup \ldots \cup \De^{(k-1)}_1)
\stackrel{\cong}\lra
(\partial \De^{(k)}_2 \backslash M_2) \cap (\De^{(1)}_2 \cup \ldots \cup \De^{(k-1)}_2)
\]
which can obviously be extended to a homeomorphism $g \colon \De^{(k)}_1 \backslash M_1 \stackrel{\cong}\lra \De^{(k)}_2 \backslash M_2$ 
such that $g (\eta_1 \backslash M_1) = \eta_2 \backslash M_2$. The list $\De^{(m+1)}_2, \ldots, \De^{(l)}_2$ must exhaust all wedges 
of $\TT_2$ because otherwise $Q_{\TT_2}$ will have more edges than $Q_{\TT_1}$ (caps do not give rise to edges of exchange 
quivers).

Part (3): We first extend the homeomorphism $g$ from the previous stage to the caps of $(S_1, M_1)$. 
Denote by $\{ \al^{(1)}_1, \ldots, \al^{(p)}_1 \}$ the arcs of $\TT_1$ that are on the boundary of exactly one 
triangle of the subsurface of $S_1$ in the lhs of \eqref{g-inter}. The map $g$ restricts to homeomorphisms
\begin{equation}
\label{g-al}
g \colon \al^{(k)}_1 \backslash M_1 \stackrel{\cong}\lra \al^{(k)}_2 \backslash M_2
\end{equation}
for some arcs $\{ \al^{(1)}_1, \ldots, \al^{(p)}_1 \}$ of $\TT_2$ -- they must be exactly those arcs of $\TT_2$ that are on the boundary of exactly one 
triangle of the subsurface of $S_2$ in the rhs of \eqref{g-inter}. Since each arc of a triangulation is on the boundary of exactly two 
triangles counted with multiplicities, $\{ \al^{(1)}_i, \ldots, \al^{(p)}_i\}$ are the boundary arcs of all caps of $\TT_i$, $i=1,2$. 
If $C^{(k)}_i$ are the caps of $\TT_i$ attached to $\al^{(k)}_i$, then the restrictions of $g$ from \eqref{g-al} can be extended 
to homeomorphisms $g \colon C^{(k)}_1 \backslash M_1 \stackrel{\cong}\lra C^{(k)}_2 \backslash M_2$ that satisify 
$g(C^{(k)}_1 \backslash M_1) \cap \partial S_1 \stackrel{\cong}\lra (C^{(k)}_1 \backslash M_1) \cap \partial S_2$.

 At this point we have a homeomorphism 
 \[
 g \colon S_1 \backslash M_1 \stackrel{\cong}\lra S_2 \backslash M_2
 \]
that satisfies \eqref{ma} and $g(\partial S_1 \backslash M_1) = S_2 \backslash M_2$. It is an easy topological fact 
that such a $g$ extends to a homeomorphism $S_1 \stackrel{\cong}\lra S_2$ that is a bijection between the sets of punctures
of $S_1$ and $S_2$ and the sets of boundary marked points of $S_1$ and $S_2$.
\end{proof} 
\noindent
{\em{Proof of \thref{homeo}.}} Denote by $P_i$ the set of punctures of $(S_i, M_i)$. 
Under the identification of $\Zset_2^{P_2}$ with the set of subsets of $P_2$, the group addition corresponds to the 
symmetric difference of subsets of $P_2$. Since $R_i$ are the subsets of punctures 
of $(S_i, M_i)$ having index $-1$ with respect to $\TT_i$, 
\[
\TT_i = (\tau(\TT_i^\circ))^{R_i}.
\]
So, the strong cluster algebra isomorphisms 
\[
\phi^\circ, \phi \in \Iso^+(\AA(S_1, M_1), \AA(S_2, M_2))
\] 
are related by 
\[
\phi = \psi_{R_2} \, \phi^\circ \, \psi_{R_1}
\]
where the automorphisms $\psi_{R_i} \in \Aut^+ \AA(S_i, M_i)$ are defined using the 
triangulations $\TT_i^\circ$, cf. \S \ref{2b.2}. By \leref{aux}, there exists an orientation-preserving 
homeomorphism $g \colon (S_1, M_1) \stackrel{\cong}\lra (S_2, M_2)$ such that $\phi^\circ = \psi_g$.
In terms of it, $\phi$ is given by
\[
\phi = \psi_{R_2} \psi_g \psi_{R_1} = \psi_{R_2} \psi_{g(R_1)} \psi_g = \psi_{R_2 \ominus g(R_1)} \psi_g = 
\psi_{g, R_2 \ominus g(R_1)}.
\]
\qed
\subsection{Uniqueness of block decompositions for maximal triangulations}
\label{4.2}
\begin{figure}[ht]
\begin{center}
\begin{tabular}{cccccc}
\begin{tikzpicture}[line cap=round,line join=round,>=triangle 45,x=1.0cm,y=1.0cm,scale=.6]
\begin{scope}[thick, every node/.style={sloped,allow upside down}]
\clip(-.5,-1) rectangle (2.5,3);
\draw(0,0) node{\Large{$\circ$}};
\draw(2,0) node{\Large{$\circ$}};
\draw(.3,0) -- node{\midarrow} (1.7,0);
\end{scope}
\end{tikzpicture}
&
\begin{tikzpicture}[line cap=round,line join=round,>=triangle 45,x=1.0cm,y=1.0cm,scale=.6]
\begin{scope}[thick, every node/.style={sloped,allow upside down}]
\clip(-.5,-1) rectangle (2.5,3);
\draw(0,0) node{\Large{$\circ$}};
\draw(2,0) node{\Large{$\circ$}};
\draw(1,1) node{\Large{$\circ$}};
\draw(1.7,0) -- node{\midarrow} (.3,0);
\draw(.2,.2) -- node{\midarrow} (.8,.8);
\draw(1.2,.8) -- node{\midarrow} (1.8,.2);
\end{scope}
\end{tikzpicture}
&
\begin{tikzpicture}[line cap=round,line join=round,>=triangle 45,x=1.0cm,y=1.0cm,scale=.6]
\begin{scope}[thick, every node/.style={sloped,allow upside down}]
\clip(-.5,-1) rectangle (2.5,3);
\draw(0,0) node{\Large{$\bullet$}};
\draw(2,0) node{\Large{$\bullet$}};
\draw(1,1) node{\Large{$\circ$}};
\draw(0,0) -- node{\midarrow} (.8,.8);
\draw(2,0) -- node{\midarrow} (1.2,.8);
\end{scope}
\end{tikzpicture}
&
\begin{tikzpicture}[line cap=round,line join=round,>=triangle 45,x=1.0cm,y=1.0cm,scale=.6]
\begin{scope}[thick, every node/.style={sloped,allow upside down}]
\clip(-.5,-1) rectangle (2.5,3);
\draw(0,0) node{\Large{$\bullet$}};
\draw(2,0) node{\Large{$\bullet$}};
\draw(1,1) node{\Large{$\circ$}};
\draw(.8,.8) -- node{\midarrow} (0,0);
\draw(1.2,.8) -- node{\midarrow} (2,0);
\end{scope}
\end{tikzpicture}
&
\begin{tikzpicture}[line cap=round,line join=round,>=triangle 45,x=1.0cm,y=1.0cm,scale=.6]
\begin{scope}[thick, every node/.style={sloped,allow upside down}]
\clip(-.5,-1) rectangle (2.5,3);
\draw(1,0) node{\Large{$\bullet$}};
\draw(1,2)node{\Large{$\bullet$}};
\draw(0,1) node{\Large{$\circ$}};
\draw(2,1) node{\Large{$\circ$}};
\draw(1,0) -- node{\midarrow} (1.8,.8);
\draw(0.2,.8) -- node{\midarrow} (1,0);
\draw(1.7,1) -- node{\midarrow} (.3,1);
\draw(.2,1.2) -- node{\midarrow} (1,2);
\draw(1,2) -- node{\midarrow} (1.8,1.2);
\end{scope}
\end{tikzpicture}
&
\begin{tikzpicture}[line cap=round,line join=round,>=triangle 45,x=1.0cm,y=1.0cm,scale=.6]
\begin{scope}[thick, every node/.style={sloped,allow upside down}]
\clip(-.5,-1) rectangle (2.5,3);
\draw(0,0) node{\Large{$\bullet$}};
\draw(0,2) node{\Large{$\bullet$}};
\draw(2,0) node{\Large{$\bullet$}};
\draw(2,2) node{\Large{$\bullet$}};
\draw(1,1) node{\Large{$\circ$}};
\draw(0,0) -- node{\midarrow} (0,2);
\draw(0,0) -- node{\midarrow} (2,0);
\draw(2,2) -- node{\midarrow} (0,2);
\draw(2,2) -- node{\midarrow} (2,0);
\draw(.8,.8) -- node{\midarrow} (0,0);
\draw(1.2,1.2) -- node{\midarrow} (2,2);
\draw(0,2) -- node{\midarrow} (.8,1.2);
\draw(2,0) -- node{\midarrow} (1.2,.8);
\end{scope}
\end{tikzpicture}\\
I&II&IIIa&IIIb&IV&V\\
\end{tabular}
\end{center}
\caption{Blocks for exchange quivers for surfaces
\label{block}
}
\end{figure}
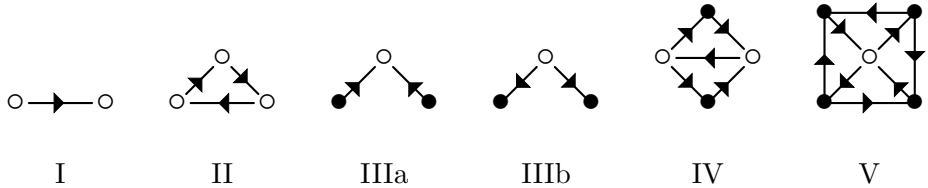

Following Fomin, Shapiro and Thurston \cite{FST}, define a {\em{block}} to be one of the six quivers on \figref{block},
called blocks of types I-V, respectively.
The unfilled vertices of the blocks are called outlets. The blocks come from puzzle pieces, which are glued groups of triangles. 
A type I block comes from a wedge, type II from a face, and type III form the puzzle piece in \figref{self-fold} when exactly one of the outer arcs
is a boundary arc. Types IV and V come from the puzzle pieces in \figref{self-fold} (a) and (b), containing self-folded triangles, when 
all outer arcs are internal. 

A quiver is called {\em{block-decomposable}} if it is connected and can be obtained by the following gluing procedure:
\begin{itemize}
\item
Consider a set of blocks and a partial matching of the total set of their outlets, where matching 
of an outlet to itself or to an outlet of the same block is not allowed. 
\item
The corresponding vertices for the partial matching are identified and a maximal set of 2-cycles is removed from the resulting
graph.
\end{itemize}
\bth{F-S-T} \cite{FST} {\em{(}}a{\em{)}} A quiver is the exchange quiver of a triangulation of a bordered 
oriented surface with marked points if and only if it is block-decomposable.

{\em{(}}b{\em{)}} Assume that $\TT_1$ and $\TT_2$ are triangulations of bordered oriented surfaces with marked points $(S_1,M_1)$ 
and $(S_2, M_2)$ such that $Q_{\TT_1} \cong Q_{\TT_2}$. If $Q_{\TT_1}$ has a unique block decomposition, then 
$(S_1, M_1)$ and $(S_2, M_2)$ are homeomorphic.
\eth

More precisely, Fomin, Shapiro and Thurston \cite{FST} define a puzzle piece to be one of the following: a face, a wedge, the 
groups of glued triangles in \figref{self-fold} (a) and (b) when at most one of the outer arcs are on the boundary of the surface $(S,M)$.
Let $\TT$ be a triangulation $(S,M)$ that is different from the triangulation in \figref{self-fold} (c) of the 4-punctured sphere and 
the triangulation in \figref{self-fold} (b) of the once-punctured digon (the latter is the only case when the triangulation in \figref{self-fold} (b) can have both of its 
outer arcs on the boundary of $S$). Then, \cite{FST}, the triangulation $\TT$ can be uniquely decomposed into a union of puzzle pieces
and caps. This decomposition leads to a block decomposition of $\TT$. In the opposite direction, every block decomposition of a quiver, leads to 
a triangulation of a bordered oriented surface with marked points $(S,M)$, which is a disjoint union of puzzle pieces and caps. 

\thref{homeo} is equivalent to the following result. The precise relationship is described at the end of the subsection.

\bpr{uniq} Assume that $(S,M)$ is a bordered surface with marked points which is different from the 
4-punctured sphere, the twice-punctured digon and the annulus with $(2,2)$ marked points on the boundary. 
If $\TT$ is a connected maximal triangulation of $(S,M)$ with at least 2 faces, 
then its exchange quiver $Q_\TT$ has a unique block decomposition which only involves blocks of types I and II.
\epr
\prref{uniq} implicitly excludes the twice-punctured monogon and  the once-punctured 2-, 3- and 4-gons
 because those surfaces do not admit maximal triangulations.
The proof of \thref{homeo} can be phrased entirely in terms of block decompositions, proving the above result instead.
We  wrote the proof of \thref{homeo} in terms of triangulations
 in order to better explain the nature of the homeomorphism $g \colon (S_1, M_1) \stackrel{\cong}\lra (S_2, M_2)$ 
 constructed from the graph isomorphism $\phi \colon Q_{\TT_1} \stackrel{\cong}\lra Q_{\TT_2}$.
 
The general uniqueness result for block decompositions stated in Theorem D is derived from the proposition in 
\S \ref{7.2}. 
 
 The exchange quivers of the connected maximal triangulation with at least 2 faces of the three surfaces, excluded 
 in \prref{uniq}, do not have unique block decompositions. For the 4-punctured sphere and the twice-punctured digon 
 the needed examples are produced in \figref{big}. The next example deals with the third surface.
 
\begin{figure}[ht]
\begin{center}
\begin{tabular}{ccc}

\begin{tikzpicture}[line cap=round,line join=round,>=triangle 45,x=1.0cm,y=1.0cm,scale=2]
\begin{scope}[thick, every node/.style={sloped,allow upside down}]
\clip(-.5,-.3) rectangle (2.5,1.5);

\draw[color=blue, line join=round, decorate, decoration={zigzag,
    segment length=4,
    amplitude=.9,post=lineto,
    post length=2pt
}](0,0)--(2,0);
\draw(1,1) node {\Large{$\bullet$}};
\draw(1,0) node {\Large{$\bullet$}};
\draw(0,0) node {\Large{$\bullet$}};
\draw(0,1) node {\Large{$\bullet$}};
\draw(2,1) node {\Large{$\bullet$}};
\draw(2,0) node {\Large{$\bullet$}};
\draw[color=blue, line join=round, decorate, decoration={zigzag,
    segment length=4,
    amplitude=.9,post=lineto,
    post length=2pt
}](0,1)--(2,1);

\draw(0,0) -- (0,1);
\draw(2,0) -- (2,1);
\draw(0,0) -- (2,1);
\draw(0,1) edge[bend right=20] (2,1);
\draw(0,0) edge[bend left=20] (2,0);

\begin{scriptsize}
\draw(.1,.5) node {$\al$};
\draw(.85,.5) node {$\be$};
\draw(2.1,.5) node {$\al$};
\draw(.5,.75) node {$\ga$};
\draw(1.5,.25) node {$\de$};
\end{scriptsize}

\end{scope}
\end{tikzpicture}
&
\begin{tikzpicture}[scale=.7]
\draw [color=blue, line join=round, decorate, decoration={zigzag,
    segment length=4,
    amplitude=.9,post=lineto,
    post length=2pt
}] (0,2) circle (2cm);
\draw (0,3) .. controls (1.5,0) and (-1.5,0) .. (0,3);
\draw (0,0) .. controls (1,1) and (1,2) .. (0,3);
\draw (0,0) .. controls (-1,1) and (-1,2) .. (0,3);
\draw (0,0) node {\Large{$\bullet$}}; 
\draw (0,3) node {\Large{$\bullet$}}; 
\begin{scriptsize}
\draw (-1,1) node {$\alpha$};
\draw (1,1) node {$\beta$};
\draw (-.2, 1.5) node {$\delta$};
\draw (0, .45) node {$\gamma$};
\end{scriptsize}
\draw (0,1) -- (0,3);
\draw (0,1) node {\Large{$\bullet$}}; 
\draw (2,0) node {};
\draw (-2,0) node {};
\end{tikzpicture} 
&
\begin{tikzpicture}[line cap=round,line join=round,>=triangle 45,x=1.0cm,y=1.0cm,scale=2]
\begin{scope}[thick, every node/.style={sloped,allow upside down}]
\clip(-.5,-.3) rectangle (1.5,1.5);

\draw(0,1) node {\Large{$\bullet$}};
\draw(1,1) node {\Large{$\bullet$}};
\draw(1,0) node {\Large{$\bullet$}};
\draw(0,0) node {\Large{$\bullet$}};

\draw(0,1) -- node{\midarrow} (1,1);
\draw(0,1) -- node{\midarrow} (0,0);
\draw(1,1) -- node{\midarrow} (1,0);
\draw(0,0) -- node{\midarrow} (1,0);
\draw(.95,0) -- node{\midarrow} (0,.95);
\draw(1,.05) -- node{\midarrow} (.05,1);

\begin{scriptsize}
\draw(0,1.2) node {$[\al]$};
\draw(1,-.2) node {$[\be]$};
\draw(0,-.2) node {$[\ga]$};
\draw(1,1.2) node {$[\de]$};
\end{scriptsize}
\end{scope}
\end{tikzpicture}\\
(a)&(b)&(c)

\end{tabular}
\end{center}
\caption{Non-uniqueness of block decompositions for $\wt{A}(2,2)$
\label{uniq-ex}}
\end{figure}
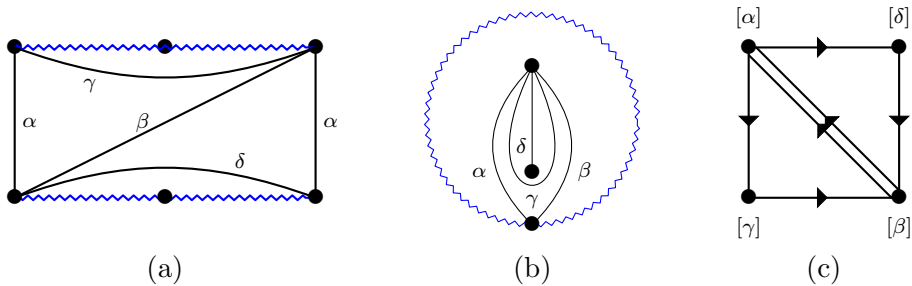

 \bex{uniq-ex} Consider the annulus with $(2,2)$ marked points on the boundary. It represents the cluster algebra 
 of type $\wt{A}(2,2)$. A maximal connected triangulation with at least 2 faces is shown in \figref{uniq-ex} (a). 
 Its exchange quiver is the one in \figref{uniq-ex} (c). It has two different block decompositions: one 
 obtained by gluing two blocks of type I, and the other by gluying a block of type I and a block of type IV. The second block decomposition 
 is associated to the triangulation of the twice-punctured monogon shown in \figref{uniq-ex} (b).
 \eex
Comparing \thref{homeo} and \prref{uniq}, we note that the annulus with $(2,2)$ marked points on the boundary is not explicitly excluded in \thref{homeo} 
because of the exclusion of the twice-punctured monogon.
\begin{proof}[Proof of \prref{uniq}] Assume that $\TT$ is a connected maximal triangulation of $(S,M)$ with at least 2 faces
whose exchange quiver $Q_\TT$ does not have a unique block decomposition. Then we can fund another surface $(S_2,M_2)$ 
and a triangulation $\TT_2$ of it such that $Q_\TT \cong Q_{\TT_2}$. The surface $(S_2,M_2)$ has to be be different from 
the twice-punctured monogon. Otherwise, an analogous argument to the one in the proof of \prref{match2} gives that 
$(S,M)$ should be isomorphic to the twice-punctured monogon or the annulus with $(2,2)$ marked points on the boundary.
Applying \thref{homeo} yields that $(S,M)$ and $(S_2, M_2)$ are homeomorphic via a homeomorphism that seds $\TT$ to $\TT_2$. 
Therefore, $\TT$ and $\TT_2$ give the same block decomposition of $Q_\TT \cong Q_{\TT_2}$, which is a 
contradiction.
\end{proof}
\subsection{A construction of connected maximal triangulations with at least 2 faces}
\label{4.3}
\bth{exist} Let $(S,M)$ be a bordered surface with marked points which is different from 
\begin{enumerate}
\item the unpunctured 4-, 5-, 6- and 7-gons, 
\item the once-punctured 2-, 3- and 4-gons,
\item the twice-punctured 1- and 2-gons, and
\item the annuli with $(1,1)$, $(2,1)$ and $(3,1)$ marked points on the boundary. 
\end{enumerate}
Then $(S,M)$ has a connected maximal triangulation with at least 2 faces.
\eth
It follows from \eqref{rank}, \eqref{fwc1}-\eqref{fwc2} and \eqref{caps} that 
\[
f(\TT) = 4(g-1) + 2b + 2 p + c(\TT) \leq 4(g-1) + 2b + 2 p + c(S,M).
\]
With the exception of the second surfaces in (2) and (3), the surfaces listed in (1)--(4) in \thref{exist} are precisely the surfaces whose 
maximal triangulations have at most one face (i.e., $4(g-1) + 2b + 2 p + c(S,M) \leq 1$).

Our notion of connected maximal triangulation with at least 2 faces is related to the notions 
of (skew-)gentle triangulations of Geiss, Labradini-Fragoso, and Schr\"oer \cite{GLS} but these notions are different. The closest
relation is between the former and the notion of gentle triangulation of \cite{GLS}, which, by definition, is one whose quiver $Q_\TT$ satisfies 
\begin{itemize}
\item its block decomposition only involves blocks of type I and II,
\item no two blocks are glued at more than one vertex,
\item all 3 cycles in $Q_\TT$ come from a single block.
\end{itemize}
An analog of our connectedness condition is not required. The proof of \thref{exist} is similar to the constriction theorem
of gentle triangulations in \cite[Theorem 7.8]{GLS} -- \thref{exist} is proved by recursively applying the following lemma
(similar to \cite[Lemma 7.9]{GLS}) starting from several initial cases. 
The difference between \thref{exist} and \cite[Theorem 7.8]{GLS} is that we have fewer initial cases.
Moreover, we only exclude finitely many surfaces while \cite[Theorem 7.8]{GLS} excludes some infinite families
and a larger (finite) set of surfaces.

\ble{add} Let $(S,M)$ be a bordered surface with marked points has a connected maximal triangulation with at least 2 faces. Assume
that $(S',M')$ is surface obtained from $(S,M)$ by adding one of the following:

{\em{(}}a{\em{)}} an additional puncture;

{\em{(}}b{\em{)}} an additional boundary component with exactly one puncture on it;

{\em{(}}c{\em{)}} an additional marked point to exactly one boundary component of $(S,M)$.
\\
Then $(S',M')$ has a connected maximal triangulation with at least 2 faces.
\ele

\begin{figure}[ht]
\begin{center}
\begin{tabular}{cccc}
\begin{tikzpicture}[line cap=round,line join=round,>=triangle 45,x=1.0cm,y=1.0cm,scale=1.5]
\begin{scope}[thick, every node/.style={sloped,allow upside down}]
\clip(-.1,-.1) rectangle (2.1,2.1);

\draw(0,0) node {\Large{$\bullet$}};
\draw(1,1) node {\Large{$\bullet$}};
\draw(1,2) node {\Large{$\bullet$}};
\draw(2,0) node {\Large{$\bullet$}};
\draw(0,0) -- (1,1);
\draw(0,0) -- (1,2);
\draw(1,2) -- (1,1);
\draw(2,0) -- (1,1);
\draw(2,0) -- (1,2);
\draw(0,0) -- (2,0);

\end{scope}
\end{tikzpicture}
&
\begin{tikzpicture}[line cap=round,line join=round,>=triangle 45,x=1.0cm,y=1.0cm,scale=1.5]
\begin{scope}[thick, every node/.style={sloped,allow upside down}]
\clip(-.1,-.1) rectangle (2.1,2.1);

\draw(0,0) node {\Large{$\bullet$}};
\draw(1,1) node {\Large{$\bullet$}};
\draw(1,2) node {\Large{$\bullet$}};
\draw(2,0) node {\Large{$\bullet$}};
\draw (1,.8) [color=blue, line join=round, decorate, decoration={zigzag,
    segment length=4,
    amplitude=.9,post=lineto,
    post length=2pt
}] circle (.2cm);
\draw(0,0) -- (1,2);
\draw(1,2) -- (1,1);
\draw(2,0) -- (1,2);
\draw(0,0) -- (2,0);
\draw (0,0) .. controls (.7,1.2) .. (1,1);
\draw (2,0) .. controls (1.3,1.2) .. (1,1);
\draw (0,0) .. controls (1.7,.7) and (1.55,1.1) .. (1,1);

\end{scope}
\end{tikzpicture}
&
\begin{tikzpicture}[line cap=round,line join=round,>=triangle 45,x=1.0cm,y=1.0cm,scale=1.5]
\begin{scope}[thick, every node/.style={sloped,allow upside down}]
\clip(-.1,-.1) rectangle (1.3,2.1);

\draw(0,0) node {\Large{$\bullet$}};
\draw(.5,1) node {\Large{$\bullet$}};
\draw(1,2) node {\Large{$\bullet$}};
\draw(.25,.5) node {\Large{$\bullet$}};

\draw(0,0) edge[color=blue, line join=round, decorate, decoration={zigzag,
    segment length=4,
    amplitude=.9,post=lineto,
    post length=2pt
}] (1,2);
\draw(0,0) edge[bend right=60] (1,2);
\draw(.25,.5) edge[bend right=45] (1,2);

\end{scope}
\end{tikzpicture}
&
\begin{tikzpicture}[line cap=round,line join=round,>=triangle 45,x=1.0cm,y=1.0cm,scale=1.5]
\begin{scope}[thick, every node/.style={sloped,allow upside down}]
\clip(-.3,-.1) rectangle (2.3,2.1);

\draw(0,0) node {\Large{$\bullet$}};
\draw(.5,1) node {\Large{$\bullet$}};
\draw(1,2) node {\Large{$\bullet$}};
\draw(2,0) node {\Large{$\bullet$}};
\draw(0,0) edge[color=blue, line join=round, decorate, decoration={zigzag,
    segment length=4,
    amplitude=.9,post=lineto,
    post length=2pt
}] (1,2);
\draw(2,0) -- (1,2);
\draw(0,0) -- (2,0);
\draw(0,0) edge[bend right=30] (1,2);
\draw(2.1,.2) node {$D$};

\end{scope}
\end{tikzpicture}
\\
(a)&(b)&(c1)&(c2)\\
\end{tabular}
\end{center}
\caption{Inductive construction of connected maximal triangulations with at least two faces. 
\label{surface build}}
\end{figure}
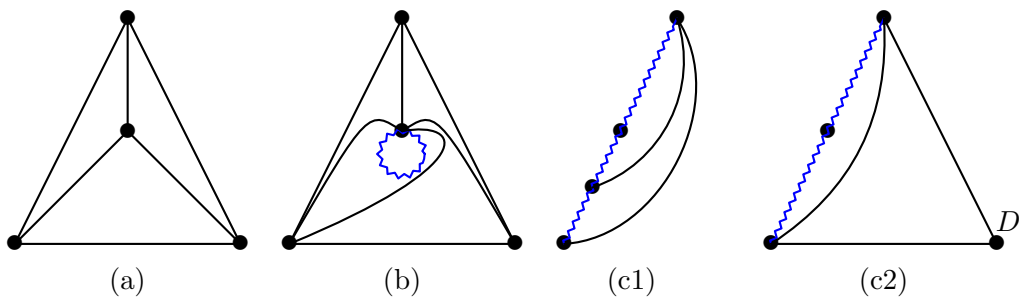

\begin{proof} Let $\TT$ be a connected maximal triangulation of $(S,M)$ with at least 2 faces.

 (a) Choose any face $\De$ in $\TT$. Place the additional puncture of $(S',M')$ in the interior of $\De$. Let $\TT'$ 
 be triangulation of $(S',M')$ obtained from $\TT$ by adding the barycentric subdivition of $\De$ by the new puncture, 
 as in \figref{surface build} (a). It is obvious that $\TT'$ a connected maximal triangulation of $(S',M')$.
 
 (b) Fix, again, a  face $\De$ of $\TT$. Let $(S',M')$ be the surface obtained from $(S,M)$ by inserting the new boundary component 
 in the interior of $\TT'$. Let $\TT'$ be the triangulation of $(S',M')$ in which $\De$ is replaced by (b) in \figref{surface build}. 
 Thus, $\TT'$ is obtained from $\TT$ by replacing the face $\De$ with three faces and a wedge.
 It is clear that $\TT'$ has no self-folded triangles and negatively double-glued pairs of triangles, and that it
 has maximal number of cups. The three new faces of $\TT'$ ensure that each of the vertices of $Q_{\TT'}$ corresponding to the 
 three new arcs of $\TT'$ are connected (via paths) to each of the vertices of $Q_{\TT'}$ corresponding to the boundary arcs of $\De'$. 
 This implies that $\TT'$ is connected.

 (c) Let $B_i \subset \partial S$ be the boundary component of $(S,M)$ in which we are adding an additional marked point. 
 Let $h_i$ be the number of marked points on this component before adding the additional point. This leads to two cases:

(c1) If $h_i$ is even, then only caps of $\TT$ are attached to $B$. We place the new marked point of $(S',M')$ in one of those caps. Let $\TT'$ 
be the triangulation of $(S,M)$ obtained from $\TT$ by subdividing the cap 
into the union of a cap and a wedge by adding an extra arc, as in \figref{surface build} (c1). 
The two boundary arcs of the wedge are obviously distinct, so $\TT'$ does not have any pairs of negatively double-glued triangles. 
The rest of the stated properties of $\TT'$ are satisfied since its faces are precisely the faces of $\TT$.
 
(c2) If $h_i$ is odd, then one wedge is attached to the boundary component $B$. Let $(S',M')$ be the surface obtained 
from $(S,M)$ by placing the new marked point on the boundary of the wedge. Denote by $D$ the vertex of the wedge that is opposite 
to the boundary arc of the wedge. Let $\TT'$ be triangulation obtained from $\TT$ by adding an arc that subdivides the 
wedge into a cap and a face, as in \figref{surface build} (c2). It is clear that the new triangulation is maximal. 
We need to prove that it is connected. Assume not; then the wedge in $\TT$ is adjacent to 2 triangles $\De_1$ and $\De_2$ 
of $\TT$ that are caps or wedges. 

If $\De_1$ or $\De_2$ is a cap, then $D \in B_i$ and the condition $c(\TT) = c(S,M)$ 
implies that both $\De_1$ and $\De_2$ are caps since $\TT$ can have at most one wedge attached to $B_i$. This implies 
that $(S,M)$ is the unpunctured hexagon which does not have a connected maximal triangulation with at least 2 faces.
This is a contradiction. 

If $\De_1$ and $\De_2$ are both wedges, then the two boundary arcs of each of them, $\eta_1$ and $\eta_2$ contain $D$. This means 
that $D$ is a marked point on a boundary component, say $B_j$, and that $B_j$ contain $\eta_1$ and $\eta_2$. 
This implies that there are 2 wedges of $\TT$ attached to $B_j$, which contradicts with the maximality condition 
 $c(\TT) = c(S,M)$.
\end{proof}

\begin{figure}[ht]
\begin{center}
\begin{tabular}{cc}
\begin{tikzpicture}[line cap=round,line join=round,>=triangle 45,x=1.0cm,y=1.0cm,scale=.35]
\clip(-8.19,-4.3) rectangle (5.6,7.66);
\draw  (-3.36,-3.56)-- (0.42,-3.58);
\draw (0.42,-3.58)[dotted]-- (3.33,-1.16);
\draw (0.42,-3.58)--(-4.97,5.88);
\draw (3.33,-1.16)--(-4.97,5.88);
\draw  (3.33,-1.16)-- (4,2.56);
\draw (4,2.56)-- (2.13,5.84);
\draw  (2.13,5.84)-- (-1.42,7.15);
\draw  (-1.42,7.15)-- (-4.97,5.88);
\draw  (-4.97,5.88)-- (-6.88,2.61);
\draw  (-6.88,2.61)-- (-6.25,-1.11);
\draw (-6.25,-1.11)-- (-3.36,-3.56);
\draw (-4.97,5.88)-- (2.13,5.84);
\draw (2.13,5.84)-- (3.33,-1.16);
\draw (0.42,-3.58)-- (-6.25,-1.11);
\draw (-6.25,-1.11)-- (-4.97,5.88);
\begin{scriptsize}
\fill[color=black] (-3.36,-3.56) node {\large{$\bullet$}};
\draw[color=black](-5.1,-2.6) node {$\alpha_g$};
\fill[color=black] (0.42,-3.58) node {\large{$\bullet$}};
\draw[color=black](-1.2,-4) node {$\beta_{g-1}$};
\fill [color=black] (3.33,-1.16) node {\large{$\bullet$}};
\fill[color=black] (4,2.56) node {\large{$\bullet$}};
\draw[color=black](4.2,.8) node {$\alpha_1$};
\fill [color=black] (2.13,5.84) node {\large{$\bullet$}};
\draw[color=black] (3.6,4.1) node {$\beta_1$};
\fill [color=black] (-1.42,7.15) node {\large{$\bullet$}};
\draw[color=black] (0.35,6.8) node {$\alpha_1$};
\fill [color=black] (-4.97,5.88) node {\large{$\bullet$}};
\draw[color=black] (-3.3,6.9) node {$\beta_g$};
\fill [color=black] (-6.88,2.61) node {\large{$\bullet$}};
\draw[color=black] (-6.4,4.3) node {$\alpha_g$};
\fill [color=black] (-6.25,-1.11) node {\large{$\bullet$}};
\draw[color=black] (-7,.9) node {$\beta_g$};
\end{scriptsize}
\end{tikzpicture}
&
\begin{tikzpicture}[line cap=round,line join=round,>=triangle 45,x=1.0cm,y=1.0cm,scale=.35]
\clip(-8.19,-4.3) rectangle (5.6,7.66);
\draw  (-3.36,-3.56)-- (0.42,-3.58);
\draw (0.42,-3.58)[dotted]-- (3.33,-1.16);
\draw (0.42,-3.58)--(-4.97,5.88);
\draw (3.33,-1.16)--(-4.97,5.88);
\draw  (3.33,-1.16)-- (4,2.56);
\draw (4,2.56)-- (2.13,5.84);
\draw  (2.13,5.84)-- (-1.42,7.15);
\draw  (-1.42,7.15)-- (-4.97,5.88);
\draw  (-4.97,5.88)-- (-6.88,2.61);
\draw  (-6.88,2.61)-- (-6.25,-1.11);
\draw (-6.25,-1.11)-- (-3.36,-3.56);
\draw (-4.97,5.88)-- (2.13,5.84);
\draw (2.13,5.84)-- (3.33,-1.16);
\draw (0.42,-3.58)-- (-6.25,-1.11);
\draw (-6.25,-1.11)-- (-4.97,5.88);
\draw (2.13,5.84)[color=blue, line join=round, decorate, decoration={zigzag,
    segment length=4,
    amplitude=.9,post=lineto,
    post length=2pt
}] .. controls (-4,4) and (2,0).. (2.13,5.84);
\draw (2.13,5.84) .. controls (2.8,0)          ..(-4.97,5.88);

\begin{scriptsize}
\fill[color=black] (-3.36,-3.56) node {\large{$\bullet$}};
\draw[color=black](-5.1,-2.6) node {$\alpha_g$};
\fill[color=black] (0.42,-3.58) node {\large{$\bullet$}};
\draw[color=black](-1.2,-4) node {$\beta_{g-1}$};
\fill [color=black] (3.33,-1.16) node {\large{$\bullet$}};
\fill[color=black] (4,2.56) node {\large{$\bullet$}};
\draw[color=black](4.2,.8) node {$\alpha_1$};
\fill [color=black] (2.13,5.84) node {\large{$\bullet$}};
\draw[color=black] (3.6,4.1) node {$\beta_1$};
\fill [color=black] (-1.42,7.15) node {\large{$\bullet$}};
\draw[color=black] (0.35,6.8) node {$\alpha_1$};
\fill [color=black] (-4.97,5.88) node {\large{$\bullet$}};
\draw[color=black] (-3.3,6.9) node {$\beta_g$};
\fill [color=black] (-6.88,2.61) node {\large{$\bullet$}};
\draw[color=black] (-6.4,4.3) node {$\alpha_g$};
\fill [color=black] (-6.25,-1.11) node {\large{$\bullet$}};
\draw[color=black] (-7,.9) node {$\beta_g$};
\end{scriptsize}
\end{tikzpicture}\\
(a)&(b)\\
\end{tabular}
\end{center}
\caption{Maximal connected triangulation with at least two faces of the once-punctured surface of genus $g$ and the unpunctured 
genus $g$ surface with one boundary component and one marked point
\label{max triang of arb genus}}
\end{figure}
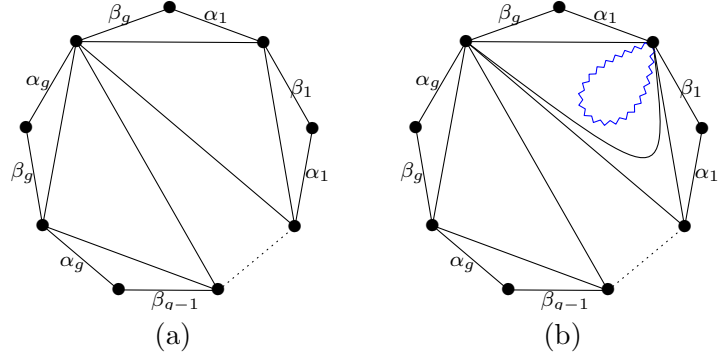

\begin{proof}[Proof of \thref{exist}] The desired triangulations of the surfaces in the theorem are obtained from maximal connected triangulations 
with at least 2 faces of either the once-punctured closed surface of genus $g \geq 1$, the unpunctured surface of genus $g \geq 1$ with one boundary component and one marked point, or one of the following 8 exceptions by recursively applying the procedures in (a)-(c) of \leref{add}:

\begin{enumerate}
\item the unpunctured 8-gon;
\item the once-punctured 5-gon;
\item the twice-punctured 3-gon;
\item the thrice-punctured 1-gon;
\item the 5 punctured sphere;
\item the unpunctured annuli with (4,1) and (2,2) marked points on the boundary;
\item the once punctured annulus with (1,1) marked points on the boundary;
\item the genus 0 surface with 3 boundary components and one marked point on each component.
\end{enumerate}
The desired triangulations for the surfaces in the exceptional list are easily constructed and are left to the reader. 
The triangulations for the two families of surfaces of arbitrary genus are shown in \figref{max triang of arb genus}.
\end{proof}

\sectionnew{Isomorphisms between surface cluster algebras and uniqueness of block decompositions}
\label{cl-iso}
In this section we give proofs of Theorems B and D from the introduction using the reconstruction theorem for maximal triangulations.
\subsection{The isomorphism problem}
\label{7.1}
\bth{iso} {\em{(}}Gu, Fomin--Shapiro--Thurston{\em{)}} \cite{G1,FST} Let $(S_1,M_1)$ and $(S_2, M_2)$ be two bordered surfaces with marked points. Then, 
$\Iso^+(\AA(S_1, M_1), \AA(S_2, M_2) \neq \varnothing$ if and only if $\Iso(\AA(S_1, M_1), \AA(S_2, M_2))  \neq \varnothing$.

Furthermore, $\Iso^+(\AA(S_1, M_1), \AA(S_2, M_2)) \neq \varnothing$ if an only if 
there exists a homeomorphism $g \colon (S_1, M_1) \stackrel{\cong}\lra (S_2, M_2)$
or the pair of surfaces is one of the following:

{\em{(}}a{\em{)}} the unpunctured hexagon and once-punctured triangle;

{\em{(}}b{\em{)}} the twice-punctured monogon and the annulus with $(2,2)$ marked points on the boundary.
\eth
\begin{proof} The equivalence 
\[
\Iso^+(\AA(S_1, M_1), \AA(S_2, M_2)) \neq \varnothing \Leftrightarrow \Iso(\AA(S_1, M_1), \AA(S_2, M_2))  \neq \varnothing
\]
follows from the existence of the non-strong automorphism $\psi_\iota \in \Aut \AA(S_2, M_2) $ from \eqref{nonstr-auto}.
For simplicity of the arguments, we will call two cluster algebras $\AA_1$ and $\AA_2$ isomorphic if 
$\Iso (\AA_1, \AA_2) \neq \varnothing$.

We group the set of all bordered surfaces with marked points into four classes:

(A) The unpuctured 4-, 5-, 6- and 7-gons, the once-punctured 2-, 3- and 4-gons,  
and the annuli with $(1,1)$, $(2,1)$ and $(3,1)$ marked points on the boundary. 

(B) The the twice-punctured monogon, and the annulus with $(2,2)$ marked points on the boundary.

(C) The twice-punctured digon and the 4-punctured sphere.

(D) All other surfaces.

Theorems \ref{thomeo} and \ref{texist} imply that the cluster algebras for the surfaces in (D) are non-isomorphic to each other. 

The remaining part of the proof is simpler. 
The cluster algebras of the two surfaces in (B) are isomorphic to each other. They are not isomorphic to any of
the cluster algebras for the surfaces in (D) -- this follows by applying \thref{homeo} to the 
annulus with $(2,2)$ marked points on the boundary and any of the surfaces in group (D).

The cluster algebras in (C) are not isomorphic to any of those in the class (D), 
because the exchange quivers with maximal number of edges of the cluster algebras in (D) have a unique 
block decomposition due to \prref{uniq}. The exchange quivers of the cluster algebras in (C) do not 
have unique block decomposition. This follows from the existence of the two triangulations in \figref{big}. 
In the case of the 4-punctured sphere we let $\eta'_1 = \eta''_1$ and $\eta'_2 = \eta''_2$. In the case 
of the twice-punctured digon we take all $\eta$'s to be boundary arcs. 

The cluster algebras of the two surfaces in (C) are not isomorphic to each other and to the cluster algebra in (B). 
The cluster algebra in (B) and the one of the twice-punctured monogon are of rank 4 and are not isomorphic to each 
other by a direct check of their triangulations. The cluster algebra of the 4-punctured sphere has rank 6 and 
is not isomorphic to any of them.

The cluster algebras in (A) do not have seeds whose quivers have two 3-cycles that share at most one common edge, 
while all other surface cluster 
algebras have this property. Because of this the cluster algebras in (A) are not isomorphic to any of the 
ones in (B)-(D). The cluster algebras in (A) have the following types, respectively: 
$A_1$, $A_2$, $A_3$, $A_4$, $A_1 \oplus A_1$,  $A_3$, $D_4$, $\wt{A}(1,1)$, $\wt{A}(2,1)$ and $\wt{A}(3,1)$. 
The isomorphisms between them (known from \cite{FZ2,FST}) are precisely the ones specified in group (a) of the statement of the theorem.
\end{proof} 
\subsection{Uniqueness of block decompositions}
\label{7.2}
\bth{uniq2} {\em{(}}Gu{\em{)}} \cite{G1} Assume that $(S,M)$ is a bordered surface with marked points which is different from one of the following:

The 4-punctured sphere, the twice-punctured digon, the once-punctured 2-, 3- and 4-gons,
the annulus with $(2,2)$ marked points on the boundary and the unpunctured hexagon.

Then the exchange quiver $Q_\TT$ of every triangulation $\TT$ of $(S,M)$ has a unique block decomposition.
\eth
\begin{proof} Assume that this is not true. Then the exchange quiver $Q_{\TT}$ needs to have a second block decomposition which would give rise 
to a bordered surface with marked points $(S_2, M_2)$, a triangulation $\TT_2$ of it, and a quiver isomorphism 
$\phi \colon Q_{\TT} \stackrel{\cong}\lra Q_{\TT_2}$. It follows from \thref{iso} that $(S_2, M_2)$ is homeomorphic to $(S,M)$.
We identify $(S_2, M_2)$ with $(S,M)$ via any of those homeomorphisms. \thref{Auto1} implies that there exist
$g \in \MCG(S,M)$ and a subset $R$ of the set of punctures of $(S,M)$ such that $\phi = \psi_{g, R} = \psi_R \psi_g$. 
Since $\TT_2$ and $\TT$ are ordinary triangulations, $R$ must be a subset of the set of punctures inside the self-folded triangles of 
$\TT_2$ and this set must be the image under $g$ of a subset of the set of punctures inside the self-folded triangles of 
$\TT$. This implies that the second block decomposition of $Q_{\TT}$ (coming from $\TT_2$)
is identical to the first block decomposition (coming from $\TT$). (The point here is that $\psi_R \colon Q_{\TT_2} \stackrel{\cong}\lra Q_{\TT_2}$ 
interchanges the vertices of $Q_{\TT_2}$ indexed by the radii of some of the self-folded triangles of $\TT_2$ with the vertices of $Q_{\TT_2}$ indexed by the corresponding loops. Such an automorphism $\psi_R$ of $Q_{\TT_2}$ does not change the blocks of $Q_{\TT_2}$, i.e., preserves its block decomposition.)
The fact that the two block decompositions are identical is a contradiction, which proves the theorem.
\end{proof}
\sectionnew{Classification of cluster automorphism groups}
\label{cl-auto}
In this section we prove Theorems C from the introduction using the main result on maximal triangulations. 
\subsection{Groups of strong cluster automorphisms}
\label{6.1}
\bth{Auto1} {\em{(}}Gu, Bridgeland--Smith{\em{)}} \cite{G1,BS} Let $(S,M)$ be a bordered surface with marked points which is different from 
the 4-punctured sphere, the once-punctured 2- and 4-gons, and the twice-punctured digon. The group of strong cluster automorphisms 
of the corresponding cluster algebra is given as follows:

{\em{(}}a{\em{)}} If $(S,M)$ is not a once-punctured closed surface, then
\[
\MCG_{\bowtie}(S,M) \cong \Aut^+ \AA(S,M) \quad \mbox{where} \quad 
(g, R) \mt \psi_{g,R}
\]
for $g \in \MCG(S,M)$ and a subset $R$ of the set $P$ of punctures of $(S,M)$.

{\em{(}}b{\em{)}} If $(S,M)$ is a once-punctured closed surface, then
\[
\MCG(S,M) \cong \Aut^+ \AA(S,M) \quad \mbox{where} \quad 
g \in \MCG(S,M) \mt \psi_g.
\]
\eth
Independently of Gu's result \cite{G1}, this theorem was stated as a conjecture in \cite{ASS} and proved in the special cases 
of disk with 1, 2 or 3 punctures and an annulus without punctures.
\begin{proof} By \eqref{ASS-emb1}--\eqref{ASS-emb2} we have an embedding of the groups on the left to the ones on the right. 
The key point in the theorem is to prove that these embeddings are surjective.
Theorems \ref{thomeo} and \ref{texist} imply the surjectivity 
of the embeddings for all bordered surfaces with marked points different from the 4-punctured sphere or 
the ones listed in (1)--(4) of \thref{exist}. It remains to establish the surjectivity of the embeddings in these 
special cases.  

(a) The surfaces in (1)--(3) of \thref{exist} (that are not excluded in the statement of this theorem) 
are the unpunctured 4-, 5-, 6- and 7-gon, the once-punctured 3-gon and the twice-punctured monogon,  
The corresponding 
cluster algebras are of type $A_1, A_2, A_3, A_3$, and $\wt{A}(2,2)$, respectively. The surjectivity for them is 
straightforward to verify by directly listing all tagged triangulations of the surfaces.

(b) The surfaces in (4) of \thref{exist} are annuli with $(1,1)$, $(2,1)$ and $(3,1)$ marked points on the boundary, 
corresponding to the cluster algebras of type $\wt{A}(1,1)$, $\wt{A}(2,1)$ and $\wt{A}(3,1)$. The surjectivity for them 
is again easily verified by listing all tagged triangulations.
\end{proof}

\bpr{auto-4sph} Let $(S,M)$ be the 4-punctures sphere. Then the index of the embedding in \eqref{ASS-emb2} equals 2:
\[
[\Aut^+(\AA(S,M)) : \MCG_{\bowtie}(S,M)] = 2.
\]
\epr

This discrepancy is due to the pathological problem with the two non-homeomorphic maximal triangulations (with the same exchange quivers) 
of the 4-punctured sphere in the proof of \prref{match2}.

\begin{proof} The triangulation in \figref{big} (a) with $\eta'_1 = \eta''_1$ is a maximal triangulation of $(S,M)$. 
Denote this triangulation by $\TT_1$, assigning a particular labeling of the punctures. 
By a direct verification one checks that 
all maximal triangulations $(S,M)$ with isomorphic exchange graphs are of the form $\sig(\TT_1)$ 
for $\sig \in A_4$ (the alternating group acting by permutations of the set of punctures).

Denote by $\TT_2$ the triangulation in \figref{self-fold} (c) for a particular labeling of its vertices. The triangulations 
of the same topological type that have isomorphic exchange quivers are of the form $\sig(\TT_2)$ for $\sig \in A_4$. 
The triangulations $\TT_1$ and $\TT_2$ have isomorphic exchange quivers as in \figref{4sphere}.
By \reref{4sph}, all non-maximal triangulations of $(S,M)$ have quivers with strictly fewer than 12 edges.
Thus, the set of tagged triangulations of the 4-punctured sphere is 
\[
X:= \{ \sig \left(  (\tau(\TT_1))^R \right), \sig \left(  (\tau(\TT_2))^R \right) \mid \sig \in A_4, R \subseteq P \}
\]
where $P$ denotes the set of 4 punctures. This implies that we have a bijection
\[
\Aut^+ \AA(S,M) \cong X, \quad \phi \in \Aut^+ \AA(S,M) \mt \phi(\TT_2).
\]
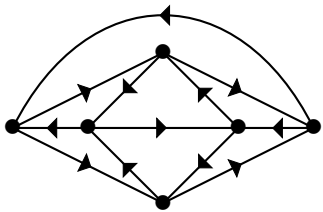
\begin{figure}
\begin{tikzpicture}[line cap=round,line join=round,>=triangle 45,x=1.0cm,y=1.0cm,scale=1]
\begin{scope}[thick, every node/.style={sloped,allow upside down}]
\clip(-.5,-1) rectangle (4.5,3);
\draw(0,1) node{\Large{$\bullet$}};
\draw(1,1) node{\Large{$\bullet$}};
\draw(2,0) node{\Large{$\bullet$}};
\draw(2,2) node{\Large{$\bullet$}};
\draw(3,1) node{\Large{$\bullet$}};
\draw(4,1) node{\Large{$\bullet$}};

\draw(0,1) -- node{\midarrow} (2,2);
\draw(0,1) -- node{\midarrow} (2,0);
\draw(1,1) -- node{\midarrow} (0,1);
\draw(1,1) -- node{\midarrow} (3,1);
\draw(2,0) -- node{\midarrow} (1,1);
\draw(2,0) -- node{\midarrow} (4,1);
\draw(2,2) -- node{\midarrow} (1,1);
\draw(2,2) -- node{\midarrow} (4,1);
\draw(3,1) -- node{\midarrow} (2,2);
\draw(3,1) -- node{\midarrow} (2,0);
\draw(4,1) -- node{\midarrow} (3,1);
\draw(4,1) .. controls (3,3) and (1,3) .. node{\midarrow} (0,1);

\end{scope}
\end{tikzpicture}
\caption{The exchange quiver of all maximal triangulations of the 4-punctured sphere
\label{4sphere}}
\end{figure}

The mapping class group of the  4-punctured sphere is $\MCG(S,M) \cong A_4$. The tagged mapping class group 
preserves the topological type of a triangulation and there is a bijection
\[
\MCG_{\bowtie}(S,M) \cong \{ \sig \left(  (\tau(\TT_2))^R \right) \mid \sig \in A_4, R \subseteq P \}, \quad 
g \in \MCG_{\bowtie}(S,M) \mt g(\TT_2).
\] 
The two bijections imply that 
\[
[\Aut^+(\AA(S,M)) : \MCG_{\bowtie}(S,M)] = 2.
\]
\end{proof}
\subsection{Groups of cluster automorphisms}
\label{6.2}
Each bordered (oriented) surface with marked points $(S,M)$ has an orientation-reversing homeomorphism 
$\iota \in \Homeo(S,M)$, $\iota \notin \Homeo^+(S,M)$. It gives rise 
to the non-strong automorphism
\begin{equation}
\label{nonstr-auto}
\psi_\iota \in \Aut \AA(S,M) , \quad \psi_\iota \notin \Aut^+ \AA(S,M) .
\end{equation} 
Composing the isomorphisms in \thref{Auto1} with $\psi_\iota$,
gives the following:
\bth{Auto2} Let $(S,M)$ be a bordered surface with marked points which is different 
from the 4-punctured sphere, once-punctured 4-gon and twice-punctured digon. 

{\em{(}}a{\em{)}} If $(S,M)$ is not a once-punctured closed surface, then
\[
\MCG_{\bowtie}^\pm(S,M) \cong \Aut \AA(S,M) \quad \mbox{where} \quad 
(g, R) \mt \psi_{g,R}
\]
for $g \in \MCG^\pm(S,M)$ and a subset $R$ of the set of punctures of $(S,M)$.

{\em{(}}b{\em{)}} If $(S,M)$ is a once-punctured closed surface, then
\[
\MCG(S,M)^\pm \cong \Aut \AA(S,M) \quad \mbox{where} \quad 
g \in \MCG^\pm(S,M) \mt \psi_g.
\]
\eth

\end{document}